
\magnification=\magstep1
\vsize=22truecm
\input amstex
\documentstyle{amsppt}
\leftheadtext{J. Jer\'onimo-Castro, E. Makai, Jr.}
\rightheadtext{Pairs of convex bodies in $S^d$, ${\Bbb R}^d$ and $H^d$}
\topmatter
\title 
\centerline{Pairs of convex bodies in $S^{\text{{\it{d}}}}$, 
${\Bbb R}^{\text{{\it{d}}}}$ 
and $H^{\text{{\it{d}}}}$,}
\centerline{with symmetric intersections}
\centerline{of their congruent copies}
\endtitle
\author J. Jer\'onimo-Castro*, E. Makai, Jr.** \endauthor
\address 
$^*$ 
Facultad de Ingenier\'\i a, Universidad Aut\'onoma de Quer\'etaro, Centro
Uni\-ver\-si\-ta\-rio, 
\newline Cerro de las Campanas s/n C.P. 76010,
Santiago de Quer\'etaro, Qro. M\'exico, ME\-XI\-CO
\newline
$^{**}$ Alfr\'ed R\'enyi Mathematical Institute, Hungarian Academy of Sciences,
\newline
H-1364 Budapest, Pf. 127, HUNGARY
\newline
{\rm{http://www.renyi.hu/\~{}makai}}
\endaddress
\email 
$^*$ jeronimo\@cimat.mx, jesusjero\@hotmail.com
\newline
$^{**}$ makai\@renyi.hu\endemail
\thanks *Research (partially) supported by CONACYT, SNI 38848
\newline
**Research (partially) supported by Hungarian National Foundation for 
Scientific Research, grant nos. T046846, T043520, K68398, K81146\endthanks
\subjclass {\it Mathematics Subject Classification 2000.} 
52A55
\endsubjclass
\abstract 
High proved the following theorem. 
If the intersections of any two congruent
copies of a plane convex body are centrally symmetric, then this body is a
circle. 
In our paper we extend the theorem of High to spherical and hyperbolic planes. 
If in any of these planes, or in ${\Bbb R}^2$, 
there is a pair of closed convex sets with interior points,
and the intersections of any congruent copies of these sets are centrally 
symmetric, then, under some mild hypotheses,
our sets are congruent circles, or, for ${\Bbb R}^2$, two parallel strips.
We prove the analogue of this statement, for $S^d$, ${\Bbb R}^d$, $H^d$,
if we suppose $C^2_+$: again, our sets are congruent balls.
In $S^2$, ${\Bbb R}^2$ and $H^2$
we investigate a variant of 
this question: supposing that 
the numbers of connected
components of the boundaries of both sets are finite, we exactly describe all 
pairs of such closed convex sets, with interior points, whose any
congruent copies have an intersection with axial symmetry
(there are 1, 5 or 9 cases, respectively).
\endabstract
\endtopmatter\document

\head 1. Introduction \endhead

By a {\it{convex body}} 
in $S^d$ (sphere), ${\Bbb R}^d$, $H^d$ (hyperbolic space)
we mean a compact convex set, with non-empty
interior. For convexity of $K \subset S^d$, with $K$ closed and 
int$\,K \ne \emptyset $,
it suffices to suppose, that for
any two non-antipodal points of $K$ the shorter great circle arc 
connecting them belongs to
$K$; then for $\pm x \in K$, $y \in $ int\,$K$ and $y \ne \pm x$,
the shorter arcs ${\widehat{(\pm x)y}}$
belong to $K$, hence some half large circle connects $\pm x$ in $K$.
In $S^d$, when saying {\it{ball}}, or {\it{sphere}}, 
we always mean one with radius at most $\pi /2$ (thus the ball is convex).
A convex body in $S^d$, ${\Bbb R}^d$, $H^d$ is {\it{strictly
convex}}, if its boundary does not contain a non-trivial segment.


R. High proved the following theorem.


\proclaim{Theorem 1} ([H])
Let $K \subset {\Bbb R}^2$ be a convex body. Then
the following statements are equivalent:

\newpage

\roster
\item
All intersections $(\varphi K) \cap (\psi K)$, 
having interior points,
where $\varphi, \psi :
{\Bbb R}^2 \to {\Bbb R}^2$ are congruences, are centrally
symmetric.
\item
$K$ is a circle. $ \blacksquare $
\endroster
\endproclaim

It seems, that his proof gives the analogous statement, when $\varphi , \psi $
are only allowed to be orientation preserving congruences.


\definition{Problem} 
Describe the pairs of closed convex sets with interior points, in $S^d$,
${\Bbb R}^d$, and $H^d$, whose any congruent copies have a centrally symmetric
intersection, provided this intersection has interior points.
Evidently, two
congruent balls (for $S^d$ of radii at most $\pi /2$), or two parallel slabs
in ${\Bbb R}^d$,
have a centrally symmetric intersection, provided it has a non-empty interior.
\enddefinition


The authors are indebted to L. Montejano (Mexico City) and G. Weiss (Dresden)
for having turned their interest to characterizations of pairs of convex
bodies with all translated/congruent copies having a centrally or axially
symmetric intersection or convex hull of the union, respectively, 
or with other symmetry properties, e.g., having some affine symmetry.


The aim of our paper will be to give partial answers to this
problem. To exclude trivialities, we always
suppose, that {\it{our sets are different from the whole plane, or space}},
and also we investigate only such cases, when the 
{\it{intersection has interior points}}.
We prove the analogue of Theorem 1 of High for $S^2$ and $H^2$.
Namely, 
we characterize the pairs of 
closed convex sets with interior points, in $S^2$, ${\Bbb R}^2$, and
$H^2$, having
centrally symmetric intersections of all congruent
copies --- provided these intersections have non-empty interiors ---   
however, for $H^2$ only under some mild hypothesis.
For $S^d$, ${\Bbb R}^d$, and $H^d$, where $d>2$,
we prove the analogous
theorem under some regularity assumptions
(weaker than $C^2$, or $C^2_+$, respectively). The only
possibilities are for $d=2$
two congruent circles, or two parallel strips for ${\Bbb
R}^2$, and for $d>2$ two congruent balls.
Moreover we investigate a variant of this question, for $S^2$, ${\Bbb R}^2$,
and $H^2$,
when we prescribe, rather than central, the 
axial symmetry of all intersections, having non-empty interiors,
but we restrict ourselves to the case 
that the numbers of connected components of the boundaries of both sets are
finite. We exactly describe all pairs of such closed convex sets with
interior points: there are $1$, $5$ or $9$ cases, respectively.

Moreover, in $S^2$, ${\Bbb R}^2$, and
$H^2$, if all small 
intersections of congruent copies of two closed convex sets with
interior points, having a non-empty interior, have some
non-trivial symmetry, then all connected components of the boundaries of the
two sets are cycles or straight lines.
For $S^d$, ${\Bbb R}^d$, and $H^d$, under the above mentioned 
regularity assumptions, if all small 
intersections of congruent copies of two closed convex sets with
interior points, having a non-empty interior, are centrally symmetric, 
then all connected components of the boundaries of the
two sets are

\newpage

congruent spheres, paraspheres, hyperspheres, or hyperplanes.
(``Small'' means here: of sufficiently small diameter.)

Surveys about characterizations of central symmetry, for convex bodies in
${\Bbb R}^d$, cf. in
[BF], \S 14, pp. 124-127, and, more recently, in [HM], \S 4.

In a paper under preparation [J-CM] 
we will give more detailed theorems about ${\Bbb
R}^d$. We will describe the pairs of closed convex sets with interior points,
whose any congruent copies have 1) a centrally symmetric intersection (provided
this intersection has interior points), without regularity hypotheses; 2)
a centrally symmetric closed convex hull of their union, also without 
regularity hypotheses.
These results will form additions to the results of the papers [So1], [So2].

\head 2. New results \endhead

We mean by a {\it{non-trivial symmetry}} a symmetry
different from the identity. Moreover, ${\text{diam}}\,(\cdot )$ 
will denote the diameter of a set.

As general hypotheses in our theorems for $d=2$ we give
$$
\cases
X {\text{ will be }}S^2, \,\,{\Bbb R}^2, {\text{ or }}H^2, \\
{\text{ and }}K,L 
\subsetneqq X {\text{ will be closed convex sets with interior points,}} \\
{\text{and }}\varphi , \psi :X \to X {\text{ will be orientation preserving
congruences. }}
\endcases
\tag *
$$

The following Theorem 2 will be the basis of our considerations for $d=2$.


\proclaim{Theorem 2}
Assume \thetag{*}. Then we have $(1) \Longrightarrow (2)$, where
\roster
\item
There exists some $\varepsilon >0$, such that 
for each $\varphi , \psi $, for 
which\,\,\,{\rm{int}}$\left( (\varphi K) \cap \right.$ 
$\left. (\psi L) 
\right) \ne \emptyset $, 
and\,\,\,{\rm{diam}}$\,\left( (\varphi K) \cap (\psi L) \right)
\le \varepsilon $, 
we have that $(\varphi K) \cap (\psi L)$
has some non-trivial symmetry.
\item
Each connected component of the boundaries of both $K$ and $L$ is a
cycle (for $X=S^2$ a circle of radius at most $\pi /2$),
or a straight line.
\endroster
In particular, if the symmetries in {\rm{(1)}} are central symmetries,
then in {\rm{(2)}} the 
connected components of the boundaries of both $K$ and $L$
are congruent.

For $X=S^2$ and $X={\Bbb R}^2$, we have $(2) \Longleftrightarrow (1)$.
Let $X=H^2$. If, both for $K$ and $L$, the infimum of the positive
curvatures of the boundary components 
is positive, and there is at most one $0$ curvature, 
then $(2) \Longleftrightarrow (1)$.
For $X=H^2$, if for, e.g., $K$, 
the infimum of the positive curvatures is $0$, or there are
two $0$ curvatures, then
we have $(2) \not\Longrightarrow (1)$. Even, we may prescribe in any way
the curvatures of the connected hypercycle or straight line
components of $K$ (with multiplicity), in the above way,
and then we can find an $L$, so that $(2)$ holds, but $(1)$ does not hold. 
\endproclaim


\proclaim{Theorem 3}
Assume \thetag{*}, and let $X=S^2$. 
Then the following statements are equivalent:

\newpage

\roster
\item
There exists some $\varepsilon >0$, such that
for each $\varphi , \psi $, for 
which\,\,\,{\rm{int}}$\left( (\varphi K) \cap \right. $
$\left. (\psi L) 
\right) \ne \emptyset $, 
and\,\,\,{\rm{diam}}$\,\left( (\varphi K) \cap (\psi L) \right)
\le \varepsilon $, 
we have that $(\varphi K) \cap (\psi L)$
has some non-trivial symmetry.
\item
There exists some $\varepsilon >0$, such that
for each $\varphi , \psi $, for 
which\,\,\,{\rm{int}}$\left( (\varphi K) \cap \right.$
$\left. (\psi L) \right) \ne
\emptyset $, 
and\,\,\,{\rm{diam}}$\,\left( (\varphi K) \cap (\psi L) \right)
\le \varepsilon $, 
we have that $(\varphi K) \cap (\psi L)$
has an axis of symmetry.
\item
$K$ and $L$ are two circles, 
of radii at most $\pi /2$.
\endroster
In particular, 
if the symmetries in {\rm{(1)}} are central symmetries, then
in {\rm{(3)}} the two circles are congruent.
\endproclaim


\proclaim{Theorem 4}
Assume \thetag{*}, and let $X={\Bbb R}^2$. 
Then the following statements are equivalent:
\roster
\item
For each $\varphi , \psi $, for 
which\,\,\,{\rm{int}}$\left( (\varphi K) \cap (\psi L) \right) \ne
\emptyset $, we have that $(\varphi K) \cap (\psi L)$
has some non-trivial symmetry.
\item
$K$ and $L$ are two circles, or one of them is a circle and the other one is a
parallel strip or a 
half-plane, or they are two parallel strips, or they are two
half-planes.
\endroster
In particular, if the symmetries in {\rm{(1)}} are central symmetries,
then in {\rm{(2)}} we have either two congruent circles, or two
parallel strips.
If the symmetries in {\rm{(1)}} are axial symmetries,
then in {\rm{(2)}}, for the case of two parallel strips, these strips 
are congruent.
\endproclaim


The following two theorems give two different characterizations for $H^2$, 
under different additional hypotheses.


\proclaim{Theorem 5}
Assume \thetag{*}, and let $X=H^2$. 
If all connected components of the boundaries of both of $K$ and $L$ are
straight lines, let their numbers be finite.  
Then
the following statements are equivalent:
\roster
\item
For each $\varphi , \psi $, for 
which\,\,\,{\rm{int}}$\left( (\varphi K) \cap (\psi L) \right) \ne
\emptyset $, we have that $(\varphi K) \cap (\psi L)$
is centrally symmetric.
\item
$K$ and $L$ are two congruent circles.
\endroster
\endproclaim


In the following theorem, the base line of a straight line is meant to be
itself.

\proclaim{Theorem 6}
Assume \thetag{*}, and let $X=H^2$. 
Then we have $(3) \Longrightarrow (2) \Longrightarrow (1)$. Supposing that
all connected components of the boundaries of both of $K$ and $L$ are
paracycles,
hypercycles or straight lines, let their total number be finite. Then
we have also $(1) \Longrightarrow (3)$. Here:
\roster
\item
For each $\varphi , \psi $, for 
which\,\,\,{\rm{int}}$\left( (\varphi K) \cap (\psi L) \right) \ne
\emptyset $, we have that $(\varphi K) \cap (\psi L)$
has some non-trivial symmetry.

\newpage

\item
For each $\varphi , \psi $, for 
which\,\,\,{\rm{int}}$\left( (\varphi K) \cap (\psi L) \right) \ne
\emptyset $, we have that $(\varphi K) \cap (\psi L)$
is axially symmetric.
\item 
We have either (A), or (B), or (C), or (D), or (E), where 
\newline (A): Any of \,$K$ and $L$ is a circle, 
a paracycle, a convex domain bounded by
a hypercycle, 
or a half-plane
--- however, if one of \,$K$ and
$L$ is a convex set bounded by
a hypercycle or is a half-plane, 
then the other one is either a circle, or a congruent copy of
the first one.
\newline 
(B): One of \,$K$ and $L$ is a circle, and the other one is
bounded either by two hypercycles, whose base lines coincide, 
or by a hypercycle, and its base line.
\newline
(C): One of \,$K$ and $L$ is a circle, of radius $r$, say, and the other one is
bounded by at least two
hypercycles or straight lines (with all base lines different), 
whose mutual distances are at least
$2r$.
\newline
(D): One of \,$K$ and $L$ is a paracycle, and the other is a parallel domain
of some straight line, for some distance $l>0$.
\newline
(E): $K$ and $L$ are congruent, and both are parallel domains of some straight
lines, for some distance $l>0$.
\endroster
\endproclaim


Now we turn to the case $d>2$.
As general hypotheses in our statements for $d>2$, we give
$$
\cases
X {\text{ will be }}S^d, \,\,{\Bbb R}^d, 
{\text{ or }}H^d, \\
{\text{ and }}K,L 
\subsetneqq X {\text{ will be closed convex sets with interior points,}} \\
{\text{and }}\varphi , \psi :X \to X {\text{ will be orientation preserving
congruences. }}
\endcases
\tag **
$$
Further, we will need
$$
\cases
{\text{Let, for each }} x \in {\text{bd}}\,K, {\text{ or each }} 
y \in {\text{bd}}\,L,{\text{there exist an }} \\
\varepsilon (x) >0, {\text{ or }} \varepsilon (y)>0, {\text{ such that }}
K, {\text{ or }}L {\text{ contains a ball of radius }} \\
\varepsilon (x), {\text{ or }} \varepsilon (y), {\text{ containing }}
x, {\text{ or }} y {\text{ in its boundary, respectively. }} 
\endcases
\tag ***
$$
and

\newpage

$$
\cases
{\text{Let, for each }} x \in {\text{bd}}\,K, {\text{ or each }} 
y \in {\text{bd}}\,L,{\text{there exist an }} \\
\varepsilon (x) >0, {\text{ or }} \varepsilon (y)>0, {\text{ such that 
the set of points of }} K, {\text{ or }} L, \\
{\text{ lying at a distance at most }} 
\varepsilon (x), {\text{ or }} \varepsilon (y), {\text{ from }} x, 
{\text{ or from }} y, \\
{\text{ is contained in a
ball }} B
{\text{ (for }} X=S^d, \,\,{\Bbb R}^d) \\
{\text{ or in a convex set }} B {\text{ bounded by a
hypersphere (for }}X=H^d), \\
{\text{with bd\,}}B {\text{ having sectional curvatures at most }} 
\varepsilon (x), {\text{ or }} \varepsilon (y), \\
{\text{and with bd\,}}B {\text{ containing }} x, {\text{ or }} y, 
{\text{ respectively.}} 
\endcases
\tag ****
$$
Clearly $C^2$ implies (***), and $C^2_+$ implies (****), and (***) implies
smoothness, and (****) implies strict convexity, respectively.


The following Theorem 7 will be the basis of our considerations for $d>2$.
Observe that in Theorem 7, (2), for ${\Bbb R}^d$ and $H^d$,
hyperplanes cannot occur, by (****).


\proclaim{Theorem 7}
Assume (**) and (***). For $X={\Bbb R}^d,\,\,H^d$ assume also (****).
Then the following statements are equivalent.
\roster
\item
There exists some $\varepsilon >0$, such that
for each $\varphi , \psi $, for 
which\,\,\,{\rm{int}}$\left( (\varphi K) \cap \right. $
$\left. (\psi L) \right) \ne \emptyset $, 
and\,\,\,{\rm{diam}}$\,\left( (\varphi K) \cap (\psi L) \right)
\le \varepsilon $, 
we have that $(\varphi K) \cap (\psi L)$
is centrally symmetric.
\item
The connected components of the boundaries of both $K$ and $L$ are congruent 
spheres (for $X=S^d$ of radius at most $\pi /2$),
or paraspheres, or congruent hyperspheres.
\endroster
\endproclaim


\proclaim{Theorem 8}
Assume (**) and (***). For $X={\Bbb R}^d,\,\,H^d$ assume also (****).
Then the following statements are equivalent.
\roster 
\item
For each $\varphi , \psi $, for 
which\,\,\,{\rm{int}}$\left( (\varphi K) \cap (\psi L) \right) \ne
\emptyset $, we have that $(\varphi K) \cap (\psi L)$ is centrally symmetric.
\item
$K$ and $L$ are two congruent balls, and, for $X=S^d$, 
their common radius is 
at most $\pi /2$.
\endroster
\endproclaim


\definition{Remark} Possibly Theorems 7 and 8 hold without any regularity
assumption. For ${\Bbb R}^d$ (where $d \ge 2$), 
in [J-CM] we will give a proof, without any
hypotheses, that (1) of Theorem 8 is equivalent to the following: $K,L$ are two
congruent balls, or are two parallel slabs. 
This is a word for word generalization of Theorem 4, case of central symmetry.
The methods of the proofs in this
paper, and in [J-CM], are completely different.
\enddefinition


\newpage

In the proofs of our Theorems we will use some ideas of [H].

\head 3. Preliminaries \endhead

For hyperbolic plane geometry we refer to [Ba], [Bo], [L], [P], for geometry
of hyperbolic
space we refer to [AVS], [C], and for elementary differential geometry we
refer to [St]. 

We shortly recall some of the concepts to be used later. In
$S^2$, $H^2$ there are the following (complete, connected, twice 
differentiable) curves of constant curvature (in $S^2$ meaning geodesic
curvature). In $S^2$ these are 
the circles, of radii $r \in (0, \pi /2]$, with (geodesic) curvature $\cot r
\in [0, \infty )$. In $H^2$, these are circles of radii $r \in (0, \infty )$,
with curvature coth$\,r \in (1, \infty )$, paracycles, with curvature $1$, and
hypercycles, i.e., distance lines, with distance $l$ from their base lines
(i.e., the straight lines
that connect their points at infinity), with curvature tanh$\,l \in (0,1)$,
and straight lines, with curvature $0$. 
Either in $S^2$, or in $H^2$ (and also in ${\Bbb
R}^2$, where we have circles and straight lines), 
each sort of the above curves have
different curvatures, and for one sort, with different $r$ or $l$, they also
have different curvatures.
The common name of these curves 
is, except for straight lines in ${\Bbb R}^2$ and $H^2$, 
{\it{cycles}}. In $S^2$ also a great
circle is called a {\it{cycle}}, but when speaking about straight lines, for
$S^2$ this will mean great circles.
An elementary method for the calculation
of these curvatures cf. in [V].

Sometimes we will include straight lines among the hypercycles.
Then the base line of a straight line is meant to be itself.

The space $H^d$ has two usual models, in the interior of the unit ball in
${\Bbb R}^d$, namely the collinear (Caley-Klein) model, and the conformal
(Poincar\'e) model. In analogy, we will speak about collinear and conformal
models of $S^d$ in ${\Bbb R}^d$, meaning the ones obtained by central
projection (from the centre), 
or by stereographic projection (from the north pole), to the
tangent hyperplane of $S^d$,
at the south pole, in ${\Bbb R}^{d+1}$. These exist of
course only on the open southern half-sphere, or on $S^d$ minus the north pole,
respectively. Their images are ${\Bbb R}^d$. 

A {\it{paraball}} in $H^d$ is a closed convex set bounded by a parasphere.

The congruences of $S^2$, ${\Bbb R}^2$ and $H^2$ can be given as follows.
The orientation preserving ones are rotations in $S^2$, rotations and
translations in ${\Bbb R}^2$, and rotations, ``rotations about an infinite
point'', and translations along a straight line (preserving this line) in 
$H^2$. The
orientation reversing ones are glide reflections in each of 
$S^2$, ${\Bbb R}^2$, and $H^2$. 

\newpage

If a non-empty closed convex set $K$ in
${\Bbb R}^2$ or $H^2$ admits a non-trivial translation, or a glide reflection
that is not a reflection, as a congruence to itself, 
then $K$ contains the closed convex hull of
the orbit of some point, w.r.t. the subgroup generated by this congruence.
Thus, $K$ contains a straight line. If a non-empty closed convex set $K$ in
$H^2$ admits a non-trivial rotation about an infinite point, 
then, by the analogous reasoning,
$K$ contains a paracycle. In most cases in our proofs, these
containments are impossible.

\head 4. Proofs of our theorems \endhead

In the proofs of our theorems by the 
{\it{boundary components of a set}} we will
mean the connected components of the boundary of that set.
For $x_1,x_2 \in X$, we write $x_1x_2$ for the {\it{distance of $x_1$ \!and
$x_2$}}, and $[x_1,x_2]$ for the {\it{segment with these endpoints}}
(supposing that $x_1,x_2$ are not antipodal points of $X=S^d$), and 
$\widehat{x_1x_2}$ for an {\it{arc of the boundary of a closed convex set with
interior points, with these end-points}}, which convex
set will be always specified.
  

\demo{Proof of Theorem 2} 
{\bf{1.}}
We begin with the proof of the implication $(1) \Longrightarrow (2)$.

{\bf{2.}}
We begin with showing that (1) implies that
both $K$ and $L$ are smooth. Then, this will imply, by convexity, that both of
them are $C^1$.

In fact, suppose,
e.g., that $K$ is not smooth. Let $x \in {\text{bd}}\,K$ be a point of
non-smoothness. Let $\alpha \in (0, \pi )$ 
denote the angle of the positively oriented half-tangents of $K$ at $x$. 

Let $y \in {\text{bd}}\,L$ arbitrary. Let $y'$
be a point of ${\text{bd}}\,L$ very close to $y$, that follows $y$ on bd$\,L$
in the positive sense.
Consider the {\it{shorter}}, 
i.e., {\it{counterclockwise arc}}
${\widehat{yy'}}$ 
{\it{of}} ${\text{bd}}\,L$. (If ${\text{bd}}\,L$ is homeomorphic to 
$S^1$; if
the connected  component of the boundary of $L$, containing $y$, 
is homeomorphic to ${\Bbb R}$, then there is
just one such arc. Here, and also later, 
when writing {\it{shorter arc}}, we mean the shorter
one in the first case, and the unique one in the second case.) 
This arc is almost like an arc in the Euclidean
plane. In particular, its map in the conformal model of $S^2$ or
$H^2$ is a very short arc (when $y$ is mapped to $0$ in the model), 
which therefore has a
total (geodesic) curvature almost $0$. So, for each point of 
the relative interior of this arc the angle of
the positively oriented half-tangents (in the conformal model, but then also
in $S^2$ or $H^2$) 
is very small. The same statement holds for ${\Bbb R}^2$ as well.

Let $x', x'' \in {\text{bd}}\,K$
be points very close to $x$, such that the smaller, say, counterclockwise
open arc $\widehat{x'x''}$ contains $x$. 
Furthermore,
we choose the points $x',x''$ so, that, additionally, for the ratio of the
distances we have $xx':xx''=b:c$, where $b,c \in (0, \infty )$ satisfy, that 
a Euclidean triangle $T$ 
with one angle $\pi - \alpha $ and adjacent sides $b,c$ is not isosceles.
Close to $x$ we have, that ${\text{bd}}\,K$ behaves almost like two (geodesic)
segments. In particular, for $x',x''$ close enough to $x$ 

\newpage

{\it{the circles with centre $x$ and with radii at most\,\,\,
$xx' + xx''$ intersect ${\text{\rm{bd}}}\,K$ just only in two points, and there
transversally}},
\newline
so, that $x$ is on the smaller arc of ${\text{bd}}\,K$ determined by these two
points.
In fact, for points of ${\text{bd}}\,K$ in
a certain neighbourhood of $x$ this follows by the differential
geometric behaviour of ${\text{bd}}\,K$, while for the points of
${\text{bd}}\,K$, not in a certain
neighbourhood of $x$ this follows by compactness of an anyhow long arc, 
with two end-points, of
${\text{bd}}\,K$ (that is embedded homeomorphically in $S^2$, ${\Bbb
R}^2$, or $H^2$), or by the fact that, outside a very long arc, with two
end-points, of
${\text{bd}}\,K$, the points of ${\text{bd}}\,K$ are approaching infinity (only
for ${\Bbb R}^2$ and $H^2$, and then only if the connected component of 
${\text{bd}}\,K$ in question is homeomorphic to ${\Bbb R}$). 

The analogous statement holds also for ${\text{bd}}\,L$, with centre of circle
$y$, and radius of circle at most $yy'$, for $y'$ sufficiently close to $y$.  

Choosing $x'x''=y'y$, there exist orientation preserving congruences $\varphi $
and $\psi$, such that $\varphi (x')= \psi (y')$, and
$\varphi (x'')=\psi (y)$, and 
$(\varphi K) \cap (\psi L)$ is bounded by the shorter
arcs ${\widehat{\varphi (x') \psi (x'')}}$ 
of ${\text{bd}}\,(\varphi K)$ and ${\widehat{\psi (y) \psi
(y')}}$ of ${\text{bd}}\,(\psi L)$. 
Thus, this intersection is almost like a Euclidean
triangle, with inner angle at $\varphi
(x)$ equal to $\pi - \alpha $, 
hence at the other two vertices $\varphi (x')=\psi
(y')$ and $\varphi (x'')=\psi (y)$ the angles between the positively oriented
half-tangents are at least about $\pi - \alpha $. In particular, this
intersection has a non-empty interior, hence has a non-trivial symmetry.
There are just three points
on the boundary of this intersection, where the angles of the positively
oriented half-tangents are at least about $\min \{ \pi - \alpha, \alpha \} $,  
at all other points these angles are about $0$. 

Then, the set $(\varphi K)
\cap (\psi L)$ is almost like a Euclidean triangle, that is not isosceles. So,
for a sufficiently small distance $x'x''=yy'$, this set has an arbitrarily
small diameter, and cannot have any
non-trivial symmetry. This is a contradiction, showing, that both $K$ and $L$
are $C^1$.

{\bf{3.}} Now, we begin the proof of the fact, that all boundary components
both of $K$ and $L$ are either cycles, or straight lines.

By compactness, and $C^1$,
on any compact arc of the boundary of $K$, or $L$, the
italicized statement from {\bf{2}} holds uniformly at the points of the
compact arc (i.e., $x$ lying in the compact arc), 
for the values of the radius at most some $\varepsilon > 0$. 
Let us consider two connected boundary components $K'$, or $L'$, of $K$, or
$L$, respectively. Let $K''$, or $L''$ be some compact arc of $K'$, or $L'$,
respectively, provided $K'$, or $L'$ is homeomorphic to ${\Bbb R}$ (and then
necessarily tends to infinity in both directions, for ${\Bbb R}^2$ and
$H^2$).
For $K'$ or $L'$ compact, i.e., when it is homeomorphic to $S^1$, and 
when necessarily $K'={\text{bd}}\,K$, or
$L'={\text{bd}}\,L$, we choose $K''$, or $L''$ equal to bd$\,K$, or bd$\,L$,
respectively. 

\newpage

Let $\varepsilon >0$ be sufficiently small. Then, we may assume 
that the italicized statement from {\bf{2}}
holds uniformly on $K''$, or $L''$ (i.e., for $x$ in $K''$, or
$L''$), respectively, for
all radius values at
most $\varepsilon $. Let $[x_1,x_2]$, or $[y_1,y_2]$ be a chord of $K$ or $L$,
respectively, of length $\varepsilon $, where $x_2$ follows $x_1$ on
${\text{bd}}\,K$ in the positive sense, and $y_2$ follows $y_1$ on
${\text{bd}}\,L$ in the negative sense.
Let $x_1,x_2 \in
K'$, and $y_1,y_2 \in L'$, with at least one of $x_1,x_2$ belonging to the
relative interior (w.r.t. $K'$) of $K''$,
and at least one of $y_1,y_2$ belonging to the relative interior 
(w.r.t. $L'$) of $L''$.
Let us choose $\varphi $ and $ \psi $
so, that $\varphi (x_i) = \psi (y_i)$ ($i=1,2$). 

First suppose, that not
both shorter arcs ${\widehat{x_1x_2}}$ and ${\widehat{y_1y_2}}$ 
are equal to the corresponding chord. 
Then, $(\varphi K) \cap
(\psi L)$ is bounded by the shorter arcs ${\widehat{\varphi (x_1) \varphi
(x_2)}}$ and
${\widehat{\psi (y_1) \psi (y_2)}}$. By the hypothesis about the arcs, this
intersection has a non-empty interior, hence has a non-trivial symmetry. 
Observe, that this intersection has just two points of non-smoothness, namely 
$\varphi (x_1) = \psi (y_1)$ and $\varphi (x_2) = \psi (y_2)$. Thus, any
non-trivial symmetry of $(\varphi K) \cap
(\psi L)$ is a central symmetry, with centre the midpoint of the segment
joining these two non-smooth points, or is an axial symmetry, either
with axis passing
through these two points, or with axis the perpendicular bisector of the
segment with endpoints these two non-smooth points.

Now, consider the case, that both above arcs are equal to the corresponding
chord, that has length $\varepsilon $. Then, $(\varphi K) \cap
(\psi L)$ may strictly contain this chord, thus, in particular, its diameter
may be not small. 
In this case, therefore, we will
consider, rather than this intersection, this common chord, 
as a degenerate closed convex set (i.e., 
with empty interior). Observe, that 
this common chord (in general not equal to $(\varphi K) \cap
(\psi L)$) has an arbitrarily small diameter, and has all three above
mentioned non-trivial symmetries.

In both cases, the intersection (in the first case above), 
{\it{or}}\,\,\,the chord (in the
second case above), has an arbitrarily small diameter, and 
has (at least)
one of the above mentioned non-trivial symmetries.
We will say, that the direction of the straight line joining the two points 
$\varphi (x_i) = \psi (y_i)$ (for $i=1,2$)
is {\it{vertical}}, and their perpendicular bisector is {\it{horizontal}}.

{\bf{4.}}
We begin with the case, when for some sequence $\varepsilon _n \to 0$, where
each $\varepsilon _n$ is sufficiently small, we have the following. 
Either $K'$, or $L'$ has a chord $[x_1,x_2]$, or $[y_1,y_2]$,
with $x_2$ following $x_1$ in the positive sense,
or $y_2$ following $y_1$ in the negative sense, and
with at least one endpoint in the relative interior (w.r.t. 
$K'$, or $L'$) of
$K''$, or $L''$, such that the following holds.
The chord $[x_1,x_2]$, or $[y_1,y_2]$ is of length $\varepsilon _n$,
and the smaller arc determined by this chord, either on $K'$, or on $L'$,  
is not symmetrical to the perpendicular
halving straight line of the chord (in particular, the respective
smaller arc is different from
the chord). 

\newpage

Then, for these $\varepsilon _n$'s, we have the
following. 
Let $[x_1,x_2]$, or $[y_1,y_2]$ be a chord of $K'$, or of $L'$, with at
least one endpoint in the relative interior (w.r.t. $K'$, or $L'$)
of $K''$, or $L''$, and of length
$\varepsilon _n$, with $x_2$ following $x_1$ on ${\text{bd}}\,K$ 
in the positive sense,
or $y_2$ following $y_1$ on ${\text{bd}}\,L$
in the negative sense, respectively.
Let $\varphi $ and $\psi $ be chosen so, that $\varphi (x_i)=
\psi (y_i)$ (for $i=1,2$), and $(\varphi K) \cap (\psi L)$ 
is bounded by the shorter
arcs ${\widehat{\varphi (x_1) \varphi (x_2)}}$ and 
${\widehat{\psi (y_1) \psi (y_2)}}$. 
(Observe that, since at least one of the arcs ${\widehat{x_1x_2}}$ and 
${\widehat{y_1y_2}}$ is
different from the respective chord, the case that $(\varphi K) \cap (\psi L)$
strictly contains this chord, and thus is degenerate, cannot occur.)
Then, the
intersection $(\varphi K) \cap (\psi L)$ has a non-empty interior,
and has an arbitrarily small diameter. Hence it
has some non-trivial
symmetry, which cannot be a symmetry w.r.t. the horizontal axis.
That is, this symmetry is a central symmetry, or is an axial symmetry with
respect to the vertical axis.

Observe, that both central symmetry, and axial symmetry with
respect to the vertical axis, cannot occur. Namely, then we would have also
an axial symmetry w.r.t. the horizontal axis, that has already been
excluded.
 
In the case of central symmetry 
the two (smaller) arcs ${\widehat{x_1x_2}}$ of ${\text{bd}}\,K$
and ${\widehat{y_1y_2}}$ 
of ${\text{bd}}\,L$, respectively, are congruent, with $x_1$
corresponding to $y_2$, and $x_2$ corresponding to $y_1$. 
In case of axial symmetry w.r.t. the vertical axis, once more the
above arcs are congruent, but now with $x_1$ 
corresponding to $y_1$, and $x_2$ corresponding to $y_2$. 

We will consider the one-sided curvatures, provided
they exist, of $K''$ at $x_i$, in the sense towards $x_{2-i}$,
and similarly, of $L''$ at $y_j$, in the sense towards $y_{2-j}$, where $x_i$
is in the relative interior of $K''$ (w.r.t. $K'$), and $y_j$
is in the relative interior of $L''$ (w.r.t. $L'$).
For both considered symmetries, the above considered two
one-sided curvatures exist and are equal
at the corresponding points, or they both
do not exist at the corresponding points. 

Now recall, that any of $x_1,x_2$, or $y_1,y_2$ 
could be any relative interior point of $K''$, or $L''$, respectively. 

First suppose the
case that, for all choices of $x_1,x_2,y_1,y_2$, we have central symmetry. Then
$\varphi (x_1)$ corresponds by the symmetry to $\psi (y_2)$. Recall that
$x_1,y_2$ could be any relative interior points of $K''$ and $L''$. Then, for
all relative interior points of $K''$ and $L''$, the considered one-sided
curvatures exist and are equal, or they do not exist for any points. However,
convex curves --- and surfaces --- are almost everywhere twice differentiable
(more exactly, the functions having, locally, in a suitable coordinate system,
these graphs, 
have Taylor series expansions, of second degree, with error
term $o(\| x \| ^2)$;
cf. [Sch], pp. 31-32, for ${\Bbb R}^d$, that extends to $S^d$ and $H^d$
by using the collinear models; observe, that we already know, that both
${\text{bd}}\,K$ and ${\text{bd}}\,L$ are $C^1$, that simplifies the condition
in [Sch]).

\newpage

This rules out the second case. Now, replacing $x_1,y_2$ by $x_2,y_1$, we
obtain the same for one-sided curvatures, but now in the opposite sense. 
Therefore, at all relative interior points of $K''$ and $L''$, the curvatures
exist and are equal.

Second suppose the
case that, for all choices of $x_1,x_2,y_1,y_2$, we have axial symmetry, with
respect to the vertical axis. Then
$\varphi (x_1)$ corresponds by the symmetry to $\psi (y_1)$. Now,
$x_1,y_1$ could be any relative interior points of $K''$ and $L''$. Then, with
this notational change, we repeat the arguments of the preceding paragraph, and
gain that, at all relative interior points of $K''$ and $L''$, the curvatures
exist and are equal.

As third case, there remains 
the case that, for some choice of $x_1,x_2,y_1,y_2$ we have
central symmetry, and for some other choice of these points
we have axial symmetry, with
respect to the vertical axis. Now, take into consideration, that an arc, or
$S^1$, is a
connected topological space, and thus products of arcs, or $S^1$'s, 
are connected
topological spaces as well. Clearly, the configurations of the points 
$x_1,x_2,y_1,y_2$ in $K' \times K' \times L' \times L'$ 
(with $x_2$ following $x_1$ in the positive sense, and $y_2$ following $y_1$
in the negative sense), where still
we suppose, that one of $x_1,x_2$ belongs to the relative interior of $K''$, 
and 
one of $y_1,y_2$ belongs to the relative interior of $L''$ (and, of course,
still $x_1x_2=y_1y_2= \varepsilon _n$), is a connected
topological space as well. Moreover, the set of configurations of the points 
$x_1,x_2,y_1,y_2$, for which one of the considered symmetry properties holds, 
is a closed subset. Further, the union of these two closed subsets is the
entire space of all above configurations of the points $x_1,x_2,y_1,y_2$.
By connectedness, these two closed subsets must intersect. That is, we must
have a configuration, that simultaneously possesses both the central symmetry,
and the axial symmetry w.r.t. the vertical axis. This, however, 
contradicts the second paragraph of {\bf{4}}.

So, the third case cannot occur. Therefore, we must have the first, or second
case. Both had the conclusion that, 
at all relative interior points of $K''$ and $L''$, the curvatures
exist and are equal.
In other words, both $K''$ and
$L''$ have equal constant curvatures, i.e., both are arcs of congruent cycles
(including entire compact cycles, i.e., circles),
or are segments.

Since $K''$, or $L''$ were arbitrary compact subarcs of $K'$, or $L'$, if
$K'$, or $L'$ were homeomorphic to $\Bbb R $ (and they were equal to
$K'={\text{bd}}\,K$, or $L'={\text{bd}}\,L$, if $K'$, or $L'$ 
was homeomorphic to 
$S^1$), we have that, in both cases, $K'$ and $L'$ are congruent cycles, or
are straight lines.

Recall, that at the beginning of {\bf{4}} we have considered the case, that 
the chord $[x_1,x_2]$, or $[y_1,y_2]$, respectively, 
is of length $\varepsilon _n$, and the
smaller arc determined by this chord, either on $K'$, or on $L'$, 
is not symmetrical to the perpendicular
halving straight line of the chord. 

\newpage

However, this contradicts the fact, that 
$K'$ and $L'$ are congruent cycles, or
are straight lines. Hence, we have obtained a contradiction.
Therefore, the case considered at the beginning of {\bf{4}} cannot occur.

{\bf{5.}} 
Thus, there remains the case that, for each sufficiently small 
$ \varepsilon >0$, both for $K''$ and $L''$, we have that
all smaller arcs of $K''$ and $L''$, having corresponding chords of length
$\varepsilon $, hence having arbitrarily small diameters, 
are symmetrical to the perpendicular halving straight line of the chord.
Observe, that this axis of symmetry halves the smaller arc, and is
perpendicular to it
at its midpoint.

Now, let $x',x''$ belong to the relative interior (w.r.t. 
$K'$) of $K''$. 
Then, there exist $x'=x_1, \ldots ,x_n=x''$ in the relative interior (with
respect to $K'$) of $K''$, following each other in the same sense, and
such, that the distance of $x_i$ and 
$x_{i+1}$ is less than $\varepsilon $ (for $i=1,\, \dots ,n-1$). 
Then, $x_i$ and 
$x_{i+1}$ are symmetrical to each other, with respect
to the perpendicular bisector of the chord $[x_i,x_{i+1}]$. Then $x_i$ and
$x_{i+1}$ are symmetrical also w.r.t. the
perpendicular bisector of some other
chord, for which the corresponding shorter arc $I'$
contains the closed shorter arc $I={\widehat{x_ix_{i+1}}}$ 
in its relative interior, $I'$ being only slightly larger than $I$ (and the
two arcs have the same midpoint). 
In particular, either
$K''$ has equal curvatures at $x_i,x_{i+1}$, or does not have a curvature at
these points. Hence, for $x',x''$, we have that 
$K''$ has equal curvatures at $x',x''$, or does not have a curvature at
these points. However, convex curves have a curvature at almost all of their
points ([Sch], pp. 31-32, cited in detail in {\bf{4}} of this proof). 
Hence the second alternative cannot
hold, i.e.,
$K''$ has a constant curvature at each of its relative interior points. 

Since $K''$ was any compact arc of $K'$
(and was equal to $K'={\text{bd}}\,K$, if $K'$ was homeomorphic to $S^1$),
we have that 
$K'$ is a $C^2$ curve of constant curvature, i.e., a cycle, or a straight line.
A similar conclusion holds for $L'$. 
This proves the implication $(1) \Longrightarrow (2)$, that is the first
statement of our theorem.

{\bf{6.}}
The particular case of (1), with central symmetry, follows easily. Let
$\varphi K$ and $\psi L$ touch each other, and push them slightly
towards each other. Then central symmetry of the new intersection implies
equality of the curvatures of the two originally touching boundary curves.

{\bf{7.}}
We turn to the third statement, i.e., to the investigation of the implication
$(2) \Longrightarrow (1)$. 

For $X=S^2$, $(2)$ clearly implies $(1)$.

For $X={\Bbb R}^2$, we  
take in consideration the following. The closed convex sets in ${\Bbb
R}^2$, whose boundaries are disconnected, are just the parallel strips.
Furthermore, the closed 
convex sets in ${\Bbb R}^2$, with connected boundaries, 
whose boundaries are cycles or straight lines, are just circles or
half-planes, respectively. Thus, any of $K$ and $L$ can be a circle, a
parallel strip, or a half-plane.

\newpage

If, e.g., $K$ is a circle, or both $K$ and $L$ are halfplanes,
then $(\varphi K) \cap (\psi L)$ is axially symmetric. 
If both $K$ and $L$ are
parallel strips, then $(\varphi K) \cap (\psi L)$ is centrally symmetric. If,
e.g., $K$ is a half-plane and $L$ is a parallel strip, then, if  
$(\varphi K) \cap (\psi L)$ has a non-empty interior, then it is unbounded.
Thus $(1)$ is satisfied in each case.

Let $X=H^2$.
Both $K$ and $L$ is either a circle, or a paracycle, or has boundary components
which are hypercycles or straight lines. 
The infimum of the positive curvatures of the boundary
components is of course the same infimum, taken only for the hypercycle
components (if there is one).
Let us first suppose that the infimum of the positive curvatures of the
hypercycle boundary components (if there is one) 
of both of $K$ and $L$ is positive, and both for $K$ and $L$
there is at most one $0$ curvature. 
That is, the distances, for which these hypercycles are distance lines, have an
infimum $c>0$, say, and there may be still at most one straight 
line component, both for $K$ and $L$. 

Let, e.g., $K_1$ and $K_2$ be two boundary components of
$K$, and let $x_1 \in K_1$, and $x_2 \in K_2$. Let $K'$ and $L'$ be defined,
as the non-empty closed convex sets (possibly with empty interiors), 
bounded by all the straight lines
for which the boundary components are distance lines, and by the at most one
straight line component. In particular, $K_1$ and
$K_2$ are distance lines for $K_1'$ and $K_2'$, with a non-negative distance.
Then the segment $[x_1,x_2]$ intersects both $K_1'$ and $K_2'$, at
points $x_1',x_2'$, and for the distances we have $x_1x_2 \ge x_1x_1'+x_2'x_2
\ge c$. This means that the distances of the different
boundary components both of $K$,
and of $L$, are bounded from below by $c$. The same holds vacuously for circles
and paracycles. Hence, if diam\,$[(\varphi K) \cap (\psi L)] <c$, then 
$(\varphi K) \cap (\psi L)$ is compact, and is 
bounded by portions of only one boundary
component of $\varphi K$, and of $\psi L$. 

Thus $(\varphi K) \cap (\psi L)$ is the intersection of
two sets, both being a circle, a paracycle, or a convex domain bounded by a
hypercycle, including a half-plane. 
Observe that a circle, and a paracycle are axially symmetric
w.r.t. any straight line passing through their centres. Thus, 
if both above sets are a circle or a 
paracycle, then their intersection is axially symmetric. There
remain the cases when one set is a convex set bounded by a hypercycle, and the
other one is a circle, a paracycle, or a convex set bounded by a hypercycle.  
In the first case an axis of symmetry of the
intersection is a straight line passing through the centre of the circle, and
orthogonal to the base line of the hypercycle. In the second case, by
compactness of the intersection, the centre of the paracycle cannot lie at an
endpoint of the base line. Therefore an axis of symmetry of the
intersection is a straight line passing through the centre of the paracycle, 
and orthogonal to the base line of the hypercycle. 

\newpage

In the third case, again by
compactness of the intersection, the base lines of the hypercycles are not
intersecting and not parallel. Therefore, the unique straight line orthogonal
to both of them is an axis of symmetry.

Now we turn to the case when the infimum of the positive
curvatures of the boundary
components of $K$ is $0$, or there are at least two $0$ curvatures 
(and the set of these curvatures, with multiplicity, is prescribed).

We begin with an example, where both bd\,$K$ and bd\,$L$ consists of two
straight lines. 
We consider the collinear model. Let 
$l_1,l_2 \subset H^2$ be parallel (but distinct) straight 
lines, with axis of symmetry $l$. Let
$x_i,y_i \in l_i$ be points symmetric w.r.t. $l$, with all six pairwise
distances at most $\varepsilon $. Then $x_1x_2y_2y_1$ is a
symmetrical quadrangle of arbitrarily small diameter, 
and is the intersection of the convex sets $K$,
bounded by $l_1,l_2$, and $L$, bounded by
the straight lines $x_1x_2, y_1y_2$. A small generic perturbation of
this quadrangle, preserving the relations $x_i,y_i \in l_i$, will have 
no non-trivial symmetry, will preserve $K$, and will perturb $L$ to a convex
set bounded by two non-intersecting and non-parallel straight lines.

If the set (with multiplicity) of the positive 
curvatures of the connected hypercycle 
components $K_i$ of $K$ is prescribed, and has infimum $c=0$, or there are at
least two $0$ curvatures, then we make the
following modification of the above example. These hypercycles are distance
lines, for distances $c_i$. We consider a closed convex set $K'$, 
bounded in the
collinear model by at most countably infinitely many chords of the model
circle, one for each $i$,
so that with at most one exception, these chords occur in pairs having
one common endpoint (possibly the set of these pairs is empty). 
Then we replace these chords by the corresponding
distance lines, outwards from $K'$. If there are two $0$ curvatures, then
the corresponding chords should 
occur in a pair, and if $c=0$, then there should be pairs
for which both distances $c_i$ 
are arbitrarily small. In the first case just copy
the above construction. In the second case we have that the hypercycles are
arbitrarily close to their base lines, and then we have two points on both of
these hypercycles, which form a convex quadrangle, with arbitrarily small
diameter, which has generically no non-trivial symmetry. 
$ \blacksquare $
\enddemo


\demo{Proof of Theorem 3}
The implications $(3) \Longrightarrow (2) \Longrightarrow (1)$
are evident. 
Last, $(1) \Longrightarrow (3)$ follows from Theorem 2. 
The particular case with central symmetries in $(1)$ follows immediately.
$\blacksquare $
\enddemo


\demo{Proof of Theorem 4}
The implication $(2) \Longrightarrow (1)$ is evident. 

For the implication $(1) \Longrightarrow (2)$ we apply Theorem 2,  
taking in consideration the following. In {\bf{7}} of the proof of Theorem 2
we have seen that 
any of $K$ and $L$ can be a circle, a parallel strip, 
or a half-plane, and with the exception 

\newpage

of the case that, e.g., 
$K$ is a half-plane and $L$ is a parallel strip, that 
$(\varphi K) \cap (\psi L)$ has some non-trivial symmetry. 
However, the case that, e.g., $K$ is a parallel strip and $L$ is a 
half-plane, contradicts $(1)$. Thus $(1) \Longrightarrow (2)$ holds.

The two particular cases, with central, or axial
symmetries in $(1)$, follow by easy discussions.
$\blacksquare $
\enddemo


\demo{Proof of Theorem 5}
{\bf{1.}} 
The implication $(2) \Longrightarrow
(1)$ is evident, so we turn to the proof of $(1) \Longrightarrow (2)$.

{\bf{2.}}
Observe that (1) of Theorem 5 implies (1) of Theorem 2, and (1) of Theorem 2
implies, by Theorem 2, that the connected components of the boundaries both of
$K$ and $L$ are congruent cycles or straight lines.

{\bf{3.}}
From {\bf{2}} we have, 
that $K$ and $L$ are two congruent circles, two paracycles,
or all their boundary components are either
congruent hypercycles,
or straight lines. However, in the case of straight lines, their total number
is finite, by the hypothesis of the theorem.

The case, 
that $K$ and $L$ are paracycles is clearly impossible. Namely, we may choose
$\varphi $ and $\psi $ so, that $\varphi K = \psi L$, and then their 
intersection is
a paracycle. However, this has exactly one point at infinity, hence is not
centrally symmetric.

We are going to show, that also the case of (finitely many)
straight lines, and the case of hypercycles is impossible.

{\bf{4.}}
First we deal with the case, when each boundary component 
both of $K$ and $L$
are straight lines, when, by hypothesis, their total number is finite.

Now, it will be convenient to use the collinear model for $H^2$. Then, in this
model, both $K$ and $L$ will be bounded by finitely many non-intersecting
chords of the boundary circle of the model. Possibly we have chords with
common end-points. Let $K_1$, or $L_1$ be some connected component of 
${\text{bd}}\,K$, or ${\text{bd}}\,L$, respectively. We may choose $\varphi $
and $\psi $ so, that $\varphi K_1 = \psi L_1 = (\varphi K) \cap (\psi L)$, 
and this
set contains the centre of the model. Thus,
$\varphi K$ and $\psi L$ lie on opposite sides of this straight
line. Let us change
$\varphi $ and $\psi $ a bit, so that in the model
$\varphi K$ and $\psi L$ rotate a very
little bit about the centre of the model. 
We will not use new notations for the new orientation preserving congruences,
but will retain the old ones $\varphi $ and $\psi $.
Let the closure ---  taken in the model with its boundary circle --- 
of the intersection
$C$ of the 
half-circles, bounded by $\varphi K_1$, or $\psi L_1$,
and containing $\varphi K$, or $\psi L$, in their new positions, respectively,
satisfy the following. It does not contain any end-point of any chord, which
in the model represents some boundary component of $\varphi K$ or $\psi L$, 
except of
course one end-point of $\varphi K_1$, and one end-point of
$\psi L_1$. 

\newpage

This can be
attained, and implies the following. The set $C$
does not intersect any
other boundary  components of $\varphi K$, or $\psi L$, 
than those, which satisfy the following properties 1) and 2): 
\newline
1) They are in
the collinear model chords of the model circle with one common end-point
with the chords
$\varphi K_1$, or $\psi L_1$, respectively, and moreover this/these common
end-point/s lie in $C$ (i.e., is/are endpoint/s of the circular arc
corresponding to $C$). 
\newline
2) From this/these connected 
component/s of the boundaries
only a half-line is in $C$.

Then, $(\varphi K) \cap (\psi L)$ is, in
the collinear
model, a sector of the model circle, a triangle, with two sides parallel in
$H^2$, and having two finite vertices, or a quadrangle, with
opposite sides parallel in $H^2$.  
The first case gives a set having
exactly one non-smooth boundary point. If it were centrally symmetric, this
boundary point would be the center of symmetry, which is a contradiction. In
the second case 
we have a set having exactly one point at infinity, hence it is not centrally
symmetric, which is a contradiction.
In the third case, if there would be a centre of symmetry, 
that would be an inner point of our set. Then
one side and its centrally symmetric image would span straight 
lines, which are not
intersecting, and not parallel. However, any two sides of this quadrangle are
either intersecting, 
or parallel. So we have a contradiction in each of the three cases.

This ends the proof of impossibility of the case, when all (finitely many)
boundary components are straight lines.

{\bf{5.}}
There remained the case, when all connected 
components of the boundaries of  
both $K$ and $L$ are congruent hypercycles. Both for $K$ and $L$, 
there is at least one such component, since $K,L \subsetneqq H^2$. Denote by
$l$ the common value of the distance, for which these hypercycles are
distance lines of their base lines.

Now, it will be convenient to consider the collinear model. Replace the image
by $\varphi $, or by $\psi $, of each
above hypercycle, for $K$, or $L$, by its base
line, respectively. These will bound closed convex sets $K_0$ and 
$L_0$ (possibly without interior points), not containing any of the image 
hypercycles. The parallel domain of $K_0$, or $L_0$, with distance $l$,
contains $\varphi K$, or $\psi L$, respectively. However, also
these parallel domains are 
contained in $\varphi K$, or $\psi L$, 
respectively. Namely, if $x \in K_0$, and the
distance of a point $z \not\in K_0$ 
from $x$ is at most $l$, then the segment $[x,z]$
intersects some boundary component of $K_0$, say, in a point $x'$. Then, the
distance of $z$ from $x'$ is at most $l$, hence $z$ lies in $\varphi K$. 

\newpage
 
Let $K_{0,1}$, or $L_{0,1}$ denote some boundary component of $K_0$, or $L_0$,
respectively. Let us suppose, that $\varphi K$ and $\psi L$ 
are in such a position, that one
end-point of $K_{0,1}$ and one
end-point of $L_{0,1}$ coincide, their other endpoints are
different, and the interiors 
of $K_0$ or $L_0$ (if not empty), lie on the opposite 
side of $K_{0,1}$ or $L_{0,1}$, as $L_{0,1}$ or $K_{0,1}$, respectively.
(This can be attained by applying some orientation preserving congruences.)
Let $K_1$, or $L_1$ denote the boundary component of $\varphi K$, or $\psi L$, 
whose base
line is $K_{0,1}$, or $L_{0,1}$, respectively (if there are two such ones,
the one that lies on the same side of $K_{0,1}$ or $L_{0,1}$, as 
$L_{0,1}$ or $K_{0,1}$, respectively). Let us consider the
intersection $M$ of the closed convex sets bounded by $K_1$ and
$L_1$ (which evidently contain $\varphi K$, or $\psi L$, respectively). 
This is bounded by some arcs of $K_1$ and $L_1$,
having one common infinite endpoint. We have $(\varphi K) \cap (\psi L) 
\subset M$. We are
going to show, that also $M \subset (\varphi K) \cap (\psi L)$. 

It will suffice to show $M
\subset \varphi K$, 
or, in other words, that $M$ lies in the parallel domain of $K_0$
with distance $l$ (the other inclusion is proved analogously). 
The straight
line $K_{0,1}$ cuts $H^2$ in two half-planes. In the half-plane containing 
$L_{0,1}$, a point belonging to $M$ clearly belongs to $\varphi K$. 
In the other half-plane, 
a point $p$ belonging to $M$ satisfies that the distance of $p$ and $L_{0,1}$
is at most $l$, hence the distance of $p$ and $K_{0,1}$
is at most $l$, hence the distance of $p$ and $K_0$
is at most $l$, as well, as was to be shown.

Thus, we have $(\varphi K) \cap (\psi L) = M$. 
The set $M$ has just one point at infinity, which implies that it cannot be
centrally symmetric.

This ends the proof of impossibility of the case, when all
boundary components are hypercycles. 
$\blacksquare $ 
\enddemo


Before passing to the proof of Theorem 6, we introduce some terminology. If we
have a topological space, $Y$, say, then we say that some property of a point
$y \in Y$ holds {\it{generically}}, if it holds outside a nowhere dense closed
subset.

If $Y$ happens to be a connected (real) analytic manifold, and $f,g:Y \to
{\Bbb R}$ are analytic functions, then either $f$ and $g$ coincide, or else
they cannot coincide on any non-empty open subset (this is the 
{\it{principle of analytic continuation}}). Otherwise said, in the 
second case, generically, for $y \in Y$, we have $f(y) \ne g(y)$.

Observe that a finite union of nowhere dense closed subsets is itself nowhere
dense and closed. In {\bf{6}} {\bf{A}} of the proof of Theorem 6 we will have 
the following
situation. On a connected (real) analytic manifold (in fact, on $H^2$)
there are finitely many, pairwise different
analytic functions, $f_1,...,f_n:Y \to {\Bbb R}$, say.
Then generically, for $y \in Y$, 
we have that $f_1(y),...,f_n(y)$ are all different.

\newpage

Before the proof we show a trigonometrical type formula in $H^2$. It is in
a sense an analogue of the law of cosines for an angle of a triangle
in $H^2$. Namely,
the law of cosines allows us, for two circles, or radii $r, R$, and distance
of centres $c$, to determine the half central angle of the arc of the circle
of radius $r$, lying in the circle of radius $R$. We will need an analogous
formula, for a circle of radius $r$, and a hypercycle, with distance $l$ from
its base line, for the half central angle of the arc of the circle
of radius $r$, lying in the convex domain bounded by the hypercycle, when the
distance $c$ of the centre of the circle and the base line of the hypercycle is
given. We consider $c$ and $l$ as signed distances.
We admit degeneration to a straight line, i.e., {\it{we admit}} 
$l=0$: then
we choose one of the half-planes bounded by this straight line (cf. below).
This formula is
surely known, but we could not find an explicit reference. Therefore we sketch
its simple proof.

So, let us consider a hypercycle, with distance $l$ from its base
line. Moreover, let us consider a circle of radius $r$, whose centre $O$ lies
at a distance $c \ge 0$ from the base line of the hypercycle. Correspondingly,
later we will consider $l$ as a signed distance, with positive sign 
determined so that we should have
$c \ge 0$ (for $c=0$ we choose the sign some way). We want to
determine the half-angle of the arc of the circle, lying in the convex domain
bounded by the hypercycle (if the hypercycle degenerates to a straight
line, then we take the half-plane bounded by it that consists 
of the points with non-positive signed distance to this straight line).
For $l$ the signed distance, we will mean our question as the determination
of the half-angle of the arc of the circle, lying in the set given by $\{ x
\in H^2 \mid {\text{dist}}\,(x,B) \le l \} $, where dist is signed distance,
and $B$ is the base line of the hypercycle.

Clearly, the intersection is non-empty if and only if $|c-l| \le r$.
{\it{At deriving our formula (**) we assume}} $\,|c-l| \le r$.

The conformal
model shows that the circle and the hypercycle have either two common points,
or they are tangent to each other, or they are disjoint
(their images are a circle, and a circular arc or segment that cuts the model
into two connected parts).

Let $C$ be
one of the common points of the circle and the hypercycle, and let 
$A$ and $B$ be the orthogonal projections of $O$ and $C$
to the base line of the hypercycle (thus $BC=l$). We let $d:=AC$.
So we have to determine the angle $\omega = \angle COA$ (for $O$ lying on the
base-line we define $\omega $ by the evident limit procedure).

\newpage

By the law of cosines
we have 
$$ 
{\text{cosh}}\,d= {\text{cosh}}\,r \cdot {\text{cosh}}\,c - {\text{sinh}}\,r
\cdot {\text{sinh}}\,c \cdot \cos \omega \,.
$$ 
Now we calculate the angle $\alpha := \angle OAC $ (for $O$ lying on the
base-line defined as a limit). Preliminarily let us suppose $l \ne 0$, that
implies $d \ne 0$.
By the law of sines we have
$$
\sin ^2 \alpha =\sin ^2 \omega \cdot {\text{sinh}}^2 r / {\text{sinh}}^2 d
\,.
$$
Last, from the right triangle $ACB$ we have 
$$
{\text{sinh}}^2 (BC) = \sin ^2 ( \pi /2 - \alpha ) \cdot {\text{sinh}}^2 d\,.
$$
So, fixing $r,c$, and supposing $\cos \omega $ as given, we determine, by
substitutions, successively, first cosh\,$d$, then $\sin ^2 \alpha $, then
${\text{sinh}}^2 (BC)$. This last expression should equal ${\text{sinh}}^2 l$.
Solving this last equation for $\cos \omega $ (which is a quadratic equation), 
we obtain, by rearranging,  
$$
\cases
\pm {\text{sinh}}\, l = {\text{cosh}}\,r \cdot {\text{sinh}}\,c 
- {\text{sinh}}\,r \cdot {\text{cosh}}\,c \cdot \cos \omega \\
={\text{cosh}}\,r \cdot {\text{cosh}}\,c \cdot
({\text{tanh}}\,c - {\text{tanh}}\,r
\cdot \cos \omega )
\,.
\endcases
\tag *
$$

We will show that here in fact we have 
$$
{\text{sinh}}\, l = {\text{cosh}}\,r \cdot {\text{sinh}}\,c 
- {\text{sinh}}\,r \cdot {\text{cosh}}\,c \cdot \cos \omega \,.
\tag **
$$
Recall that $|c-l| \le r$ is assumed.

In \thetag{*} 
the expression in the middle lies in $[{\text{sinh}}\,(c-r),
{\text{sinh}}\,(c+r)]$. So, for $0 \le r \le c$, it is non-negative, and,
since the signed distance $l$ was taken to be positive, so that we should have
$c \ge 0$, therefore here the first
expression must be ${\text{sinh}}\,l$, i.e., we have \thetag{**}. 
Now let $r > c \ge 0$. Then the boundary
of our circle intersects the base line in two points. 
The case $l=0$ corresponds to
a well known formula for a right triangle in $H^2$: it is equivalent to
${\text{tanh}}\,c= {\text{tanh}}\,r \cdot \cos \omega $. In particular, (*)
and (**) are
valid for $l=0$ as well.
Let us increase 
$\omega $, and thus the middle expression of \thetag{*}. 
Then the signed distance of the end-point of the radius of our
circle, enclosing an angle $\omega $ with the radius of our circle
orthogonally intersecting the base line (for $O$ on the base line this is
meant as a limit), increases. This corresponds to the
fact that we have ${\text{sinh}}\,l$ in the first expression in 
\thetag{*},
i.e., we have \thetag{**}.
Last we extend the validity of \thetag{**} 
to $c < 0$. Let us apply \thetag{**} to $-c,-l, \pi - \omega
$ rather than $c,l, \omega $. Then the validity of \thetag{**} for these
values implies its validity for $c,l, \omega $. 

\newpage

Later, in the proof of Theorem 6, we will consider the case when $l \ge 0$;
then, of course, $c$ varies in ${\Bbb R}$.

Observe that \thetag{**} implies the necessary and sufficient 
condition for the existence of
a point of intersection, namely $|c-l| \le r$. In fact, the right hand side of
\thetag{**} lies in $[{\text{sinh}}\,(c-r),{\text{sinh}}\,(c+r)]$.

From this there follows the converse implication. Namely:
if \thetag{**} is satisfied, then $|c-l| \le r$, and 
there exist 1) a
hypercycle, having a signed distance $l$ from its
base line $B$, and 2) a 
circle of radius $r$, that has a centre at a distance $c$ 
from the base line of the hypercycle, such that 3) 
the circle intersects
$\{ x \in H^2 \mid {\text{dist}}\,(x,B) \le l \} $
in a circular arc of half central angle $\omega $. 

\demo{Proof of Theorem 6}

{\bf{1.}} 
The implication $(2) \Longrightarrow (1)$ is trivial.

{\bf{2.}}
We continue with the proof of $(3) \Longrightarrow (2)$.

We begin with case (A).

A circle is axially symmetric w.r.t. a straight line
spanned by any of its diameters. A paracycle is axially symmetric with respect
to any straight
line passing through its centre (its point at infinity), 
and a convex domain bounded by
a hypercycle, or a half-plane is axially symmetric w.r.t. 
any straight line, that intersects its base line, or its boundary,
orthogonally, respectively. These imply, that if any of $\varphi K$ 
and $\psi L$ is either a circle or a
paracycle, then their intersection is axially symmetric w.r.t. 
(any) straight line joining their centres. If one of $\varphi K$ and $\psi L$ 
is a circle, and the other one is
a convex set bounded by a hypercycle, or is a half-plane, then the straight
line passing
through the centre of the circle, and orthogonal to the base line of the
hypercycle, or to the boundary of the half-plane, 
is an axis of symmetry of the intersection. 

Last, let 
$\varphi K$ and $\psi L$ be congruent convex sets, both bounded by
hypercycles, or let them be two half-planes. 
Consider the base lines of these hypercycles, or the boundaries of 
these half-planes, respectively. There
are four cases. These lines
\newline 
a) may coincide; or
\newline 
b) may intersect; or
\newline 
c) may have a common point at infinity (but are distinct); or 
\newline 
d) may have no common finite or infinite point. 
\newline Case a) is evident. In case b),
$\varphi K \ne \psi L$, and
${\text{bd}}\,(\varphi K)$ and ${\text{bd}}\,(\psi L)$ intersect
transversally at some point $p$ (for this use the conformal model). 
Then, $(\varphi K) \cap (\psi L)$ has an inner
angle at $p$, of measure less than $ \pi $, and the
halving straight line of this angle is an axis of symmetry of 
$(\varphi K) \cap (\psi L)$. In case c), if one of $\varphi K$ and $\psi L$
contains the other, the intersection is evidently axially
symmetric. 

\newpage

Otherwise, the symmetry axis of the base lines is an axis of
symmetry of the intersection.
In case d), we consider the pair of points on the base lines, realizing
the distance of these lines. The straight line connecting these points
is orthogonal to both
lines, and is an axis of symmetry of $(\varphi K) \cap (\psi L)$.

We continue with case (B). If $K$ is a circle, and $\psi L$ is bounded
by two hypercycles, whose base lines coincide
(one of them possibly degenerating to a straight line), 
then the straight line passing through
the centre of $\varphi K$, and orthogonal to the above base line, is an axis of
symmetry of $(\varphi K) \cap (\psi L)$.

We continue with case (C). If $K$ is a circle of radius $r$, and
the boundary hypercycle or straight line 
components of $L$ have pairwise distances at least $2r$, then 
int\,$(\varphi K)$ can intersect at most one boundary component of $\psi L$. 

If int\,$(\varphi K)$ does not intersect any boundary component of $\psi L$
(and, by hypothesis, int\,$[(\varphi K) \cap (\psi L)] \ne \emptyset $),
then $(\varphi K) \cap (\psi L) = \varphi K$ is a circle, hence is axially
symmetric.

If int\,$(\varphi K)$ intersects exactly one boundary component $L_1$ 
of $\psi L$,
then $(\varphi K) \cap (\psi L)$ is the same as the intersection of 
$\varphi K$ and of 
the closed convex set, bounded by $L_1$, and containing $\psi
L$. This has an axis of symmetry, cf. case (A).

We continue with case (D). Let $K$ be a paracycle, and $L$ a parallel
domain of some straight line, with some distance $l>0$. 
Consider the common base line of the two
hypercycles, bounding $\varphi K$. If the infinite point of the paracycle
$\varphi K$ lies on this common base 
line, then this straight
line is an axis of symmetry of $(\varphi K) \cap (\psi L)$. 
If the infinite point of $\varphi K$ 
does not lie on this common base line, then there is a unique
straight line that passes through the infinite point of $\varphi K$, 
and is orthogonal
to the common base line. Then this unique straight line is an axis of
symmetry.
 
Last we turn to case (E). Consider the common base lines of the two
hypercycles bounding $\varphi K$, and of the two
hypercycles bounding $\psi L$. These two straight 
lines can coincide, or can intersect, or can be
parallel (but distinct), 
or can be neither intersecting nor parallel. In any case there is
an axial symmetry interchanging these two straight lines. This axial symmetry
interchanges the parallel domains of these straight lines, with distance $l$,
as well.
Hence it is a symmetry of the intersection of these parallel domains, i.e., of 
$(\varphi K) \cap (\psi L)$.

{\bf{3.}}
Last we turn to the proof of $(1) \Longrightarrow (3)$.
By Theorem 2 
we know, that each boundary component of 
both $K$ and $L$ is 
either a cycle, or a straight line.
Thus, for each of $K$ and $L$, we have the following possibilities: it is a
circle, or a paracycle, or its boundary components are 
hypercycles and straight lines. 

We make a case distinction. Either both bd\,$K$ and bd\,$L$ are connected, or
one of them has several connected components.

\newpage

{\bf{4.}}
We begin with the case when both bd\,$K$ and bd\,$L$ are connected, and will
show that then we have (A) of (3) of the theorem.

We have to investigate the cases when 
\newline 
a) $\varphi K$ and $\psi L$ are 
one paracycle and one convex
set bounded by a hypercycle or a straight line, or 
\newline 
b) $\varphi K$ and $\psi L$ are two incongruent convex sets, both
bounded by a hypercycle or a straight line,
\newline
and in both cases we have to find a contradiction.

Now it will be convenient to use the conformal model.
In case a), let the centre of the
paracycle be one endpoint of the base line of the hypercycle, or one endpoint
of the straight line. Then
$(\varphi K)\cap (\psi L)$ has a smooth boundary, except at one point $p$, 
that is
the intersection of ${\text{bd}}\,(\varphi K)$ and ${\text{bd}}\,(\psi
L)$. In case b), let the base lines of the two
hypercycles, or the base line of the hypercycle and the straight line
intersect, respectively (two straight lines cannot occur). Then, also 
${\text{bd}}\,(\varphi K)$ and ${\text{bd}}\,(\psi L)$ intersect, at a single
point, and this
point $p$ is the only non-smooth point of $(\varphi K) \cap 
(\psi L)$. 

Both in case a) and b), any non-trivial symmetry of $(\varphi K) \cap (\psi
L)$
would have $p$ as a fixed point. Thus, it would be an axial symmetry, 
w.r.t. the
angle bisector of the inner angle of $(\varphi K) \cap (\psi L)$ at $p$. Thus,
this symmetry should
interchange the portions of the boundaries of $\varphi K$ and $\psi L$, 
bounding $(\varphi K) \cap (\psi L)$. However, these portions of boundaries
have different
curvatures, which is a contradiction. 

Thus, in the case investigated in {\bf{4}}, we have shown (A) of (3) of the
theorem.

{\bf{5.}}
There remained the case when one of $K$ and $L$ 
has at least two boundary components.
Observe that this rules out the cases when
$K,L$ are two circles, or two paracycles, or one circle and one 
paracycle. There remain the cases when one of $K$ and $L$ is bounded by
hypercycles and straight lines, and the other one is a circle, of
some radius $r$, or when one of $K$ and $L$ is bounded
by finitely many
hypercycles and straight lines, and the other one is either
a paracycle, or also is bounded
by finitely many
hypercycles and straight lines. We will investigate these three cases
separately.

If one of $K$ and $L$ is bounded
by hypercycles and straight lines, then the boundary components $K_i$ of
$\varphi K$, or the boundary components $L_i$ of $\psi L$ have a natural cyclic
order, in the positive sense, on bd\,$(\varphi K)$, or bd\,$(\psi L)$,
respectively. We associate to $\varphi K$, or to $\psi L$ a graph, whose
vertices
are the infinite points of the $K_i$'s, or $L_i$'s, and between two such
points there is an edge, if they are the two infinite points of some $K_i$, or
$L_i$, respectively. We say that this edge is $K_i$, or $L_i$,
respectively. This graph can be a union of vertex-disjoint paths, or
can be a cycle. 

\newpage

Here we admit a cycle of length $2$, shortly a $2$-cycle,
when the graph consists
of two vertices, and two edges between these two vertices, which are
two $K_i$'s ($L_i$'s) with both infinite points common. 

If we have two edges in these graphs with a common vertex, and they are e.g.,
$K_1$ and $K_2$, then by this notation we will mean 
that $K_2$ follows $K_1$ on bd\,$(\varphi K)$ in the positive sense.
If $K_1$ and $K_2$ 
form a $2$-cycle, then the notation is fixed some way.
(A similar convention holds for $L$).

When one of $K$ and $L$ is bounded by hypercycles and straight lines, 
then, later in the proof, for
brevity we will write {\it{hypercycle}} for a hypercycle, or for a
straight line; i.e., {\it{the curvature is allowed to be $0$}} 
as well. The {\it{base
line}} of a straight line is considered to be the straight line itself.

{\bf{6.}}
Let $K$ be a circle of radius $r$ and centre $O$, and $L$ be bounded by
hypercycles and straight lines. We have to show that either $L$ is bounded by
two hypercycles with common base line (i.e., (B) of (3) of the theorem holds), 
or $L$ has at least two boundary
components (which holds by the assumption in the beginning of {\bf{5}}) and
these boundary components have pairwise distances at least $2r$ (i.e., (C) of
(3) of the theorem holds). Let us
suppose the contrary, i.e., that we have both that $L$ is not bounded by
two hypercycles with common base line, and that
dist\,$(L_1,L_2) <2r$, for some different boundary components 
$L_1$ and $L_2$ of $\psi L$. By dist\,$(L_1,L_2) <2r$, we have, 
for some $\varphi $, that 
int\,$(\varphi K)$ intersects both $L_1$ and $L_2$. 

In this case, int\,$(\varphi K)$ intersects $L_i$ for some $i$'s, and $\varphi
K$ touches $L_i$ for some other $i$'s. We will show that by an arbitrarily
small perturbation of the centre $O$ of $\varphi K$
we can attain that $(\varphi K) \cap
(\psi L)$ has no non-trivial symmetry. Clearly then 
we need not care those $L_i$'s, for which $(\varphi K) \cap L_i = \emptyset $.
(Observe that any compact set in $H^2$ intersects only finitely many $L_i$'s.)

{\bf{A.}}
Then $(\varphi K) \cap (\psi L)$ is 
a convex body, 
bounded, alternately, by (at least two) non-trivial arcs of bd\,$(\varphi
K)$, and (at least two) non-trivial arcs of some $L_i$'s, for different
$L_i$'s. The curvatures of these arcs are greater than $1$, or smaller than
$1$, respectively, so each congruence of $( \varphi K) \cap ( \psi L)$
preserves both above types of arcs, separately. Now consider 
conv\,$[\left( \right.$bd\,$ \left. (\varphi K) \right) 
\cap (\psi L) ]$, that is preserved
by each congruence of $( \varphi K) \cup ( \psi L)$. It is 
obtained from the circle $\varphi K$, by cutting off disjoint
circular segments by several, but at least two disjoint, non-trivial chords,  
having endpoints the points of intersection of 
the single $L_i$'s with bd\,$(\varphi K)$, 
for all $L_i$'s intersecting int\,$(\varphi K)$.
We want to attain that
$$
{\text{all these chords are of different lengths.}}
\tag 1
$$ 

\newpage

We are going to show that {\it{this can be attained by a small, generic  
motion of the centre $O$ of}} $\varphi K$.

The lengths of the above
chords of $\varphi K$ are uniquely determined by the half central angles
corresponding to the chords. We use formula \thetag{**} before the proof of
this theorem. We have several such equations, corresponding to several $L_i$'s,
with respective values
$l_i$ and $c_i$ (but with $r>0$ fixed). We have $l_i \ge 0$ for each $i$,
and then $c_i$ will vary in ${\Bbb R}$.
We will use arbitrarily small generic
perturbations of the centre $O$ of our circle $\varphi K$. Then the set of
hypercycles that intersect the perturbed
int\,$(\varphi K)$ is a subset of the set of all those
hypercycles $L_i$, that intersect a fixed concentric closed
circle (concentric meant
before perturbation) 
of some radius $r'>r$. This second set is
finite (cf. the second paragraph of {\bf{6}}); 
let it be $\{ L_i \mid i \in I \} $. 
Hence, it suffices to exclude all pairwise equalities
of finitely many expressions for $\cos \omega $ --- obtained from solving 
the equations \thetag{**}, 
before the proof of this theorem, for all $i \in I$ --- namely 
those of the form
$$
({\text{cosh}}\,r \cdot {\text{sinh}}\,c_i - {\text{{sinh}}\,l_i ) /
({\text{sinh}}\,r \cdot \text{cosh}}\,c_i) \,.
$$
Observe that all these expressions are analytic in $O$, since the $c_i$'s are
analytic in $O$ (and $r$ and the $l_i$'s are fixed). 

Moreover, none of these equations is an identity. Namely, we can consider 
a circle $\varphi K$ outside of the 
convex set bounded by the 
boundary component $L_i$ of $\psi L$, containing $\psi L$,
where $i \in I$. By a
certain motion we may attain that $\varphi K$ just touches
this $L_i$, and is otherwise outside of the convex set in the last sentence.
Then the $i$'th expression for $\cos \omega $ has value $1$, but
all other $j$'th expressions, where $j \in I$, 
have values for $\cos \omega $ not in $[-1,1]$. Hence
the $i$'th expression and any other $j$'th expression, for $i,j \in I$, 
are not identical.

Therefore all our finitely many analytic equations are not identities.
Hence each of them holds only for $O$ belonging to a 
nowhere dense closed subset. Therefore,
except for $O$ belonging to
a nowhere dense closed subset, none of our equations hold. That is, we have
proved what was claimed in \thetag{1}.

{\bf{B.}}
From now on we will suppose \thetag{1}.
There are 
two possibilities. Either we can have at least three such chords --- as in
\thetag{1} --- or we always have exactly two such chords. 

\newpage

Any congruence of $(\varphi K) \cap (\psi L)$ to itself preserves $\varphi K$, 
and conv\,$[\left( \right. $bd\,$\left. (\varphi K) \right) 
\cap (\psi L)]$, and
also the above mentioned at least three, or exactly two 
disjoint chords, since their lengths are different. However, a non-trivial
congruence preserving a single chord is an axial symmetry w.r.t. the
orthogonal bisector straight
line of the chord. However, there are no three disjoint
circular segments, cut off by chords orthogonal to a single straight line, 
containing the centre of the circle (namely, orthogonal
to their common orthogonal bisector). 

There remains the case when we can have only exactly two disjoint
circular segments cut off by chords of $\varphi K$. These must 
correspond to the
above considered hypercycles $L_1$ and $L_2$. 
Then $\varphi K$ is not even touched by
any other $L_i$, since then by a small motion of $\psi L$ we could attain that 
int\,$(\varphi K)$ intersects at least three $L_i$'s, which case was
above settled.

The above reasoning gives that, in this case, the orthogonal bisecting straight
lines
of the two chords coincide, furthermore, contain $O$; moreover,
these remain true after an arbitrary small
motion of the centre $O$ of $\varphi K$, except those into 
a nowhere dense closed subset, cf. {\bf{A}}.

However, the orthogonal bisecting straight
lines of these chords are orthogonal to the
base lines of $L_1$ and $L_2$. We have that these base lines are different,
since their coincidence was excluded in the first paragraph of {\bf{6}}. 
Then they have no finite point in common, and they have either one, or no
infinite point in common. 

If they have one infinite point in common, then they
admit no common orthogonal straight line. 

If they have no infinite points in common,
then they have exactly one common orthogonal straight
line, that should contain the
centre $O$ of $\varphi K$, for {\it{all}} \,\,small motions of $O$, 
except those into a nowhere dense closed subset, cf. {\bf{A}}. 
This is clearly impossible.

Thus we have proved, what was promised in the beginning of {\bf{6}}: 
namely that, in the case investigated in {\bf{6}},
we have cases (B) or (C) of (3) of our theorem.

{\bf{7.}}
Now let $K$ be a paracycle, and let $L$ be bounded
by finitely many, but at least two hypercycles and straight lines.

If the graph of $\psi L$ 
consists of paths, then, using one end-point of one path,
and the adjacent edge of the graph, 
we choose $\varphi K$ in the conformal
model as a circle of small radius, that is thus far from all other boundary
components of $\psi L$. Then
we can repeat the consideration from
{\bf{4}}, case a), and we obtain a contradiction.

There remains the case when the graph of $\psi L$ 
is a cycle. Hence there are two
edges $L_1$, $L_2$, with a common vertex, in the graph of $\psi L$.

\newpage

We consider the conformal model. Fixing the position of $L_1$ and $L_2$, 
we consider their common vertex (or one of their common vertices)
at infinity, $1$, say. Let $L_2$ follow $L_1$ at $1$ in the positive sense. 
We choose $\varphi K$ so that 
it touches the boundary of the model either 1) at $1$, or 2)
very close to $1$ but not at $1$, and its interior
intersects both $L_1$ and $L_2$, and
its image in the model is a circle of very small radius. Then $\varphi K$ 
is far from
all other boundary components of $\psi L$. Hence, $(\varphi K) \cap (\psi
L)$ is either 1) an arc triangle, bounded by an arc of $L_2$, an arc of
bd\,$(\varphi K)$, and an arc of $L_1$, or 2) 
an arc quadrangle, bounded by an arc of $\varphi K$, an arc of $L_2$, another 
arc of $\varphi K$, 
and an arc of $L_1$, in this cyclic order, in the positive sense. 

In case 1) $(\varphi K) \cap (\psi L)$ has exactly one infinite point, hence
its only non-trivial symmetry is an axial symmetry w.r.t. an axis passing
through this infinite point, and interchanging $L_1$ and $L_2$ (observe that
a rotation about this infinite point is impossible, by the last paragraph 
of \S 3). Hence $L_1$
and $L_2$ are congruent.
 
In case 2) neither of the 
diagonals can be an axis of symmetry, and there is no symmetry that would be
combinatorially a $4$-fold rotation. So, a non-trivial symmetry of 
$(\varphi K) \cap (\psi L)$ can be a central symmetry, or an axial symmetry
w.r.t. the common orthogonal bisector straight lines 
of the opposite arc-sides. 

We begin with the case of central symmetry. Then the opposite arc-sides of 
$(\varphi K) \cap (\psi L)$ on the paracycle $\varphi K$
are centrally symmetric images of
each other. Let us suppose that the centre of symmetry is the centre of the
(conformal) model. Then the hyperbolic central symmetry coincides with the
Euclidean central symmetry. The paracycle in the model is a Euclidean 
circle touching the boundary of
the model in one point. The centrally symmetric image of the paracycle w.r.t.
the centre of the model intersects the paracycle only in at most two points,
thus in no arc.
 
We continue with the case of axial symmetry 
w.r.t. the common orthogonal bisector straight lines 
of the opposite arc-sides, lying on
bd\,$(\varphi K)$. However, a common orthogonal bisector straight line
to the two opposite (thus disjoint)
arc-sides of $(\varphi K) \cap (\psi L)$, lying on
bd\,$(\varphi K)$, cannot
exist. Namely, such bisectors are different and
parallel, thus have no finite point in common, but have
only an infinite point in
common, namely the infinite point of $\varphi K$.

We continue with the case of axial symmetry 
w.r.t. the common orthogonal bisector straight lines 
of the opposite arc-sides, lying on
$L_1$ and $L_2$. However, a
common orthogonal straight line
to $L_1$ and $L_2$ is a common orthogonal to their base
lines as well. These base lines are parallel, hence admit no common
orthogonal, unless they coincide. If they coincide, then $\psi L$ is bounded
just by $L_1$ and $L_2$. 

The considerations in 1) and 2) yield (D) of (3) of the theorem.

\newpage

{\bf{8.}}
Last, let both $K$ and $L$ be bounded
by finitely many
hypercycles and straight lines. 
Then, by the first paragraph of {\bf{5}}, 
either $K$, or $L$ has at least two boundary components. 

We will show that the graphs of $\varphi K$ and $\psi L$ 
must have only a few edges, and we
will clarify the structure of these graphs, till we will obtain that we must
have case (E) of (3) of our theorem.

We will make the following case distinction, for the graphs of $\varphi K$
and $\psi L$. 
\newline 1) Both graphs contain a pair of edges with at least one common
end-point.
\newline 2) One graph contains a pair of 
edges with at least one common endpoint,
but the other graph does not contain such a pair, i.e., it consists of vertex
disjoint edges, and the number of these edges is at least $2$.
\newline 3) One graph contains a pair of 
edges with at least one common endpoint,
but the other graph does not contain such a pair, i.e., it consists of vertex
disjoint edges, and the number of these edges is $1$.
\newline 4) None of the graphs contains a pair of 
edges with at least one common 
endpoint, i.e., both of them
consist of vertex disjoint edges. Here, by {\bf{5}}, at
least one of the graphs contains at least $2$ edges. 
\newline These cases are exhaustive, and mutually exclusive.

{\bf{9.}}
We begin with the proof of case 1).

Then each of the graphs of $\varphi K$ and $\psi L$ 
contains a path of length
$2$ or a $2$-cycle.
The corresponding boundary components of $\varphi K$ and $\psi L$ are
denoted by $K_1,K_2$, and $L_1,L_2$, with $K_2$ following $K_1$ 
on bd\,$(\varphi K)$, and $L_2$ following $L_1$ on bd\,$(\psi L)$,
according to the positive orientation. (If one of the graphs is a $2$-cycle,
this does not determine $K_1$ etc.; then we fix some notation.)

We are going to show that $K_1 \cup K_2$ and $L_1 \cup L_2$ 
are images of each other by an orientation reversing congruence, moreover, that
the graphs of $\varphi K$ and $\psi L$ both are $2$-cycles. Then
these will imply that 
bd\,$(\varphi K)=K_1 \cup K_2$, and bd\,$(\psi L)=L_1 \cup L_2$.
  
We use the conformal model.
Recall that any three different points on the boundary of the model can be
taken by a congruence to any other three different 
boundary points of the model.
Therefore we may suppose the following. The considered
common vertex of $K_1$ and $K_2$ is
$1$, and their other vertices are very close to $-1$ --- hence all other
boundary components of $\varphi K$ are very close to $-1$, as well ---
and the considered common vertex of $L_1$ and $L_2$ is
$i$, and their other vertices are very close to $-i$ --- hence all other
boundary components of $\psi L$ are very close to $-i$, as well
(and $K_2$ follows $K_1$
on bd\,$(\varphi K)$
at $1$ in the positive sense, and similarly for $L_2,L_1$ at $i$). 

\newpage

Then we have that the distance of $0$ to 
$(\varphi K) \cap (\psi L)$ is small (possibly $0$),
and $(\varphi K) \cap (\psi L)$ is
bounded by arcs of $K_1,L_1,K_2,L_2$, in this cyclic order, in the
positive sense. In fact, all other boundary components, both of $\varphi K$, 
and of
$\psi L$, are in the model very close to the boundary of the model, hence
cannot cut off parts of this arc-quadrangle, which is
not close to the boundary.

Thus $(\varphi K) \cap (\psi L)$ is a compact arc-quadrangle. Its possible
non-trivial symmetries are, combinatorially, the following: 
two four-fold rotations, one central symmetry, and axial symmetries
w.r.t. diagonals or common orthogonal bisector straight lines 
of two opposite edges.
If we have a symmetry that is combinatorially a four-fold rotation, then we
also have a symmetry that is a combinatorial central symmetry. Hence we need
not exclude the case of a combinatorial four-fold rotation, exclusion of a 
combinatorial central symmetry will suffice.
Observe that a non-trivial symmetry of $(\varphi K) \cap (\psi L)$
extends, by analycity, 
to a non-trivial symmetry of $K_1 \cup K_2 \cup L_1 \cup L_2$.

{\bf{A.}}
We begin with the case when $(\varphi K) \cap (\psi L)$
has a central symmetry. Then $K_1$ and $K_2$ have two common infinite
points (images of each other by this symmetry), and the same statement holds
for $L_1$ and $L_2$, hence the graphs of both 
$\varphi K$ and $\psi L$ 
are $2$-cycles. Clearly then the curvatures of $K_1$ and $K_2$, as well as 
those of $L_1$ and $L_2$, are equal, and are positive.

We are going to show that also the curvatures of $K_i$ and $L_i$ coincide,
i.e., we have case (E) of (3) of our theorem.
Let, e.g., the curvature of $K_i$ be less than the curvature of $L_i$. Let us
choose a new position for $\psi L$, in such a way that the the infinite points
of each of $K_1,K_2,L_1,L_2$ are $\pm 1$, and $K_1,L_1$ are in the lower
half-plane, and $K_2,L_2$ are in the upper half-plane.
Let us rotate $\psi L$ about the infinite point $1$, 
counterclockwise, by a
rotation of a small measure. Then, in the (conformal) model, 
the Euclidean 
tangents of $L_1$ and $L_2$ at $1$ do not change during this rotation. 
Therefore, in the new position, 
$(\varphi K) \cap (\psi L)$ is an arc-triangle, bounded by
$K_2,L_2,K_1$, in this order. This has a unique infinite point, hence a
non-trivial symmetry must be an axial symmetry w.r.t. an axis passing through
this infinite point, i.e., the common base line for $K_1$ and $K_2$, i.e., the
real line. Then
this axial symmetry should preserve $L_2$, and therefore 
this axis should be the
orthogonal side
bisector of the arc-side of $(\varphi K) \cap (\psi L)$ on $L_2$. However, one
of the angles 
of this axis and $L_2$ equals one of their angles 
at their other common point $1$
in the model (with its boundary circle), that is acute.
Hence, in the new position,
$(\varphi K) \cap (\psi L)$ has no non-trivial symmetry, a contradiction. 

\newpage

That is, we have obtained case (E) of (3) of our theorem.

{\bf{B.}} 
We continue with the case when $(\varphi K) \cap (\psi L)$ has 
an axial symmetry w.r.t. a
common orthogonal bisector straight line of two opposite edges.

Observe that a common orthogonal straight line
to opposite sides, say, on $K_1$ and $K_2$, 
is a common orthogonal to
the respective base lines, which are parallel, or coinciding. In the case when
these base lines do not coincide, this is impossible. 

If these base lines
coincide, then the graph of $\varphi K$ is a $2$-cycle, and $\varphi K$ is
bounded just by $K_1$ and $K_2$. Let the other common infinite point of
$K_1$ and $K_2$ be $-1$
(and let $\psi L$ be as described in the fourth paragraph of {\bf{9}}).
The axis of our
symmetry is a straight line orthogonal to 
this common 
base line, i.e., to the real axis. 
This axis of symmetry must contain the common infinite point
$i$ of $L_1$ and $L_2$. Thus this axis of symmetry is the imaginary axis.
Now a small rotation 
of $\psi L$ about $i$ will make $\psi L$ not symmetric
w.r.t. the imaginary axis. From now on we will consider this small
perturbation of the arc-quadrangle 
$(\varphi K) \cap (\psi L)$, rather than the original arc-quadrangle.
Clearly, in the
new position, this axial symmetry becomes destroyed. 
Moreover, for all sufficiently small (depending on the measure of the above
small rotation) new perturbations of $(\varphi K) \cap (\psi
L)$, in its new position, this axial symmetry will not exist. (Else the new
position itself would have this axial symmetry.)

{\bf{C.}}
There remained the case, when the first
perturbed $(\varphi K) \cap (\psi L)$ has an 
axial symmetry w.r.t. some diagonal. Since this
symmetry of $(\varphi K) \cap (\psi L)$ yields a congruence
between $K_1 \cup K_2$ and $L_1 \cup L_2$, these last sets are also axially
symmetric images of each other, hence are images of each other by an
orientation reversing congruence, as promised.

Thus the graphs of $\varphi K$ and of $\psi L$ either 
\newline
a) simultaneously contain paths of length $2$, namely $K_1K_2$ and $L_1L_2$, 
or 
\newline
b) are simultaneously $2$-cycles, with edges $K_1,K_2$, and $L_1,L_2$,
respectively.

In case a), recall that at the beginning of {\bf{9}},
$\varphi K$ and $\psi L$ were chosen as follows: the considered
common infinite point of $K_1$ and $K_2$ is $1$, the considered
common infinite point of $L_1$ and $L_2$ is $i$, and the other end-points of
$K_1$ and $K_2$ are close to $-1$, and the other end-points of
$L_1$ and $L_2$ are close to $-i$. Observe that the small rotation applied in
{\bf{B}} preserves these properties.

Then $(\varphi K) \cap (\psi L)$ is a
compact arc-quadrangle, such that
the distance of $0$ to it is small (possibly is $0$). 
If one of the diagonals is the
axis of symmetry, then it is an angle bisector of the angles at the vertices
that it connects. Choose the arc-sides with one endpoint
at one of these vertices, say, the arc-sides on $K_2$ and $L_1$. These
determine this diagonal uniquely. Then choose a third arc-side, on $K_1$, say. 

\newpage

By symmetry, this third arc-side
already uniquely determines the fourth arc-side. However, fixing the
hypercycles containing these three arc-sides (i.e., $K_1,K_2,L_1$), 
by a small rotation of $L_2$ about $i$, that extends to an orientation
preserving congruence of $\psi L$, preserving $L_1$,
we attain that this symmetry is
destroyed. By sufficiently smallness of the second perturbation, the second
perturbed $(\varphi K) \cap (\psi L)$ does not have the symmetry investigated
in {\bf{B}}, and, by construction, does not have the symmetry investigated in 
{\bf{C}} either. This is a contradiction, hence case a) (i.e., when 
$K_1K_2$ and $L_1L_2$ are paths of length $2$) cannot exist.
(Recall that case {\bf{A}} was settled above: then we obtained case (E) of (3)
of our theorem.)

Summing up: unless we have case (E) of (3) of our theorem (in {\bf{A}}), 
we have case b) 
(i.e., that the graphs of both $\varphi K$ and $\psi L$ 
are $2$-cycles, with edges $K_1,K_2$, and $L_1,L_2$, respectively), 
and also that $K$ and $L$ are
images of each other by an orientation reversing congruence, as promised in the
beginning of {\bf{9}}.

{\bf{10.}}
We investigate further the situation described at the end of {\bf{9}}. 
By a suitable notation, 
we have that the curvatures of $K_1$ and $L_1$ are equal, and also
that the curvatures of $K_2$ and $L_2$ are equal.
We are going to show that 
$K_1$ and $K_2$, and then also $L_1$ and $L_2$, each has the same
curvature, i.e., that (E) of (3) of our theorem holds. 
Let us suppose that the curvature of $K_2$ is greater than the
curvature of $K_1$.

Observe that $\varphi K$, or $\psi L$, 
is symmetrical w.r.t. any straight line
orthogonal to the common base line of the $K_i$'s, or of the $L_i$'s,
respectively. 
Hence any two congruent copies 
of $K$ (and of $L$)
are simultaneously directly and indirectly congruent. Since $K$ and $L$ are
indirectly congruent by {\bf{9}}, they are also directly congruent.

Let us fix $\varphi K$ so that its points at infinity are $\pm 1$, and $K_1$
lies in the closed
lower half-plane, and $K_2$ lies in the closed upper half-plane (with at least
one of them lying in the respective open half-plane). Let us
obtain $\psi L$ by rotating $\varphi K$ about $1$ in positive sense a bit,
with the image of $K_i$ being $L_i$. Then $(\varphi K) \cap (\psi L)$ is
bounded by two arcs (in the model, with its boundary circle), 
one lying on $L_2$, the other one lying on
$K_1$. Then $(\varphi K) \cap (\psi L)$ has one non-smooth point, at $L_2 \cap
K_1$. Hence a non-trivial symmetry of $(\varphi K) \cap (\psi L)$
is an axial symmetry, w.r.t. the angle bisector of the angle of 
$(\varphi K) \cap (\psi L)$ at this point. However, this symmetry
interchanges the two portions of the boundary, lying on $L_2$ and
$K_1$. This contradicts the fact that the curvatures of $L_2$ and $K_1$ are
different. 

That is, we have obtained case (E) of (3) of our theorem.

{\bf{11.}}
Now we turn to the proof of case 2) from {\bf{8}}.

So let, e.g., the graph of $\varphi K$ 
consist of vertex disjoint edges, whose number is at least $2$. Let us
choose two vertex-disjoint edges of this graph, $K_1$,
$K_2$, say. 

\newpage

Further, let the graph of $\psi L$ 
contain a path of length $2$ or a $2$-cycle, 
consisting of $L_1$ and $L_2$, 
where $L_2$ follows $L_1$ in the positive orientation (if we have
a $2$-cycle, then their numeration is done in some way).
We are going to show that this case cannot occur.

We fix $\varphi K$ and thus $K_1$ and $K_2$, and will choose $\psi L$ in the
following way. The set $\varphi K$ lies in the convex set bounded by $K_1$ and
$K_2$. Then we have relatively open arcs $I_1$ and $I_2$
of the boundary of the (conformal) model, bounded by 
the infinite points of $K_1$ and $K_2$, and lying outside of 
the above mentioned convex set. We choose the (considered)
common infinite point of $L_1$
and $L_2$ at the midpoint of $I_1$, and the other infinite points of $L_1$ and
$L_2$ (possibly coinciding) very close to the centre of $I_2$.

Then $(\varphi K) \cap (\psi L)$ is contained in 
a compact arc-quadrangle $Q$, bounded by arcs
lying on $K_1,L_2$, $K_2,L_1$, 
in this order, say. 

Observe that all boundary components of $\psi L$, other than
$L_1$ and $L_2$, 
are in the model very close to the boundary of the model, hence
cannot cut off parts of this arc-quadrangle $Q$, which arc-quadrangle is
not close to the boundary. 
So these boundary components have no arcs on 
bd$\,[(\varphi K) \cap (\psi L)]$. So we need not deal with these 
boundary components.

However, there may exist several boundary components $K_i$ of $\varphi K$,
with $i\ne 1,2$,
which cut off parts of this arc-quadrangle, hence have non-trivial arcs on
bd$\,[(\varphi K) \cap (\psi L)]$. 

Since we investigate case 2), we have that the 
$K_i$'s have no common endpoints.
Of course, $L_1$ and $L_2$ have at least one common endpoint. However, by
con\-struc\-tion, 
neither $L_1$ and any $K_i$ (including $K_1$ and $K_2$),
nor $L_2$ and any $K_i$ (including $K_1$ and $K_2$), have any
common endpoint.

We are going 
to show that any non-trivial congruence of $(\varphi K) \cap (\psi L)$
is a congruence of $Q$ as well; moreover, it is a congruence preserving the
opposite pairs of arc-edges of $Q$.

We have that bd\,$[(\varphi K) \cap (\psi L)]$ consists of arcs, following
each other, in the positive sense, lying on $K_1,L_2,K_{i(1)},L_2,K_{i(2)},...,
K_{i(j)},L_2,K_2,L_1,K_{i(j+1)},L_1,K_{i(j+2)},$
\newline
$...,K_{i(k)},L_1$, say. From all of these arcs
only those lying on $L_1$ and $L_2$ lie on different hypercycles, which have
at least one point in common. 

Let us introduce a symmetric relation $\Cal R$ on the
arc-sides of $(\varphi K) \cap (\psi L)$. For two such arc-sides $S_1,S_2$ 
we have $S_1 {\Cal R} S_2$, 
if the hypercycles spanned by these sides are
different, and have at least one common end-point. Clearly 
any non-trivial congruence of $(\varphi K) \cap (\psi L)$ preserves this
relation $\Cal R$, hence also the set of arc-sides ${\Cal S} := 
\{ S_1 \mid \exists S_2
$ such that $S_1 {\Cal R} S_2 \} $. 

\newpage

Observe that the relation
$\Cal R$ induces a complete bipartite graph on the vertex set ${\Cal S}$, with
classes ${\Cal L}_i$, for 
$i=1,2$, where ${\Cal L}_i$ is the set of arc-sides
of $(\varphi K) \cap (\psi L)$, lying on $L_i$, for $i=1,2$.

Therefore, each non-trivial congruence of $(\varphi K) \cap (\psi L)$
preserves the two-element set
$\{ {\Cal L}_1, {\Cal L}_2 \} $. Of course, also
the cyclic order of the arc-sides of $(\varphi K) \cap (\psi L)$
is preserved, 
up to inversion. Let the first end-point of the first arc-side and 
last endpoint of the last arc-side in ${\Cal L}_i$ (i.e., lying on $L_i$),
be $v_{i,1}$ and $v_{i,2}$.
Then 
the set $\{ \{ v_{1,1},v_{1,2} \} , \{ v_{2,1},v_{2,2} \} \} $ is preserved by
each non-trivial congruence of $(\varphi K) \cap (\psi L)$,
as well. So, $Q$ is preserved by
each non-trivial congruence of $(\varphi K) \cap (\psi L)$
as well, even in such a way, that separately the opposite
pairs of sides are preserved.

{\bf{12}}.
By {\bf{11}} we need to discuss only the congruences of the arc-quadrangle
$Q$, more exactly only those of them, that
preserve the opposite pairs of sides. Therefore,
combinatorially, the possible non-trivial congruences, to be investigated,
and to be excluded, are central symmetry,
and axial symmetries w.r.t. common orthogonal side-bisector straight lines 
of opposite sides.

{\bf{A.}}
We begin with the case of central symmetry of $(\varphi K) \cap (\psi L)$. 
The central symmetry interchanges the arc-sides lying on $K_1$ and $K_2$,
hence also $K_1$ and $K_2$, hence also
the infinite points of $K_1$ and the infinite points of $K_2$. Thus its
centre must be the intersection of the straight lines connecting the
interchanged end-points. (This determines the interchanged pairs of
end-points uniquely.)
We may assume that this centre of symmetry is $0$.
Also, by central symmetry,
the graph of $\psi L$ must be a $2$-cycle, with $L_1$ and $L_2$ having the same
curvatures. Then the centre of symmetry must lie on
the common base line of $L_1$ and $L_2$. 

{\bf{B.}}
We continue with the case of axial symmetry of $(\varphi K) \cap (\psi L)$
w.r.t. the common orthogonal
bisector of the arc-sides lying on $L_1$ and $L_2$. Such a common orthogonal
straight
line is orthogonal to the base lines of $L_1$ and $L_2$ as well, hence it
exists only if the base lines of $L_1$ and $L_2$ coincide (they cannot be
parallel but different), i.e., the graph of $\psi L$ is a $2$-cycle. 
Then the symmetry interchanges $K_1$ and $K_2$, hence
its axis is the unique axis of symmetry interchanging $K_1$ and $K_2$. 
We may suppose that this axis is the imaginary axis. 
Hence the axis of this symmetry is orthogonal to
the common base line of $L_1$ and $L_2$. 

{\bf{C.}}
We continue with the case of axial symmetry of $(\varphi K) \cap (\psi L)$
w.r.t. the common orthogonal
bisector straight line
of the arc-sides lying on $K_1$ and $K_2$. This axis is the unique
straight line orthogonal to $K_1$ and $K_2$ (and hence also
to their base lines). 
Since the centre of symmetry
considered in {\bf{A}} was $0$, and the axis of symmetry
considered in {\bf{B}} was the imaginary axis, this axis is 
the real axis. 

\newpage

Then the axis of the unique axial symmetry of $L_1 \cup L_2$,
interchanging $L_1$ and $L_2$, is the real axis.

Considering all three possible cases {\bf{A}}, {\bf{B}}, {\bf{C}}, 
we have the following.
By a small generic perturbation of $\psi L$ we can attain that the axis of
the axial symmetry of $L_1 \cup L_2$, interchanging $L_1$ and $L_2$, 
intersects the imaginary axis at a point different
from $0$ (thus does not contain $0$), and the angle enclosed by this axis of
symmetry and the imaginary axis is different from $\pi /2$. Then we also have
that this axis of symmetry is different from the real axis. Thus this
perturbation simultaneously destroys all three possible non-trivial symmetries
of $(\varphi K) \cap (\psi L)$, discussed in {\bf{A}}, {\bf{B}}, and {\bf{C}}.
This is a contradiction. Hence, the 
case investigated in {\bf{11}} cannot occur, as promised in the first
paragraph of {\bf{11}}.

{\bf{13.}}
We turn to the proof of case 3) from {\bf{8}}.

Let, e.g., the graph of $\varphi K$ 
consist of a single edge $K_1$, i.e., this is the unique
boundary component,
and let the graph of $\psi L$ contain a
path of length $2$ or a $2$-cycle. 
We are going to show that this is impossible.

Let the graph of $\psi L$ contain two edges
$L_1$ and $L_2$ with a common vertex, 
following each other in the positive orientation, at the considered common
vertex. We consider the conformal model.
We consider $L_1$ and $L_2$ as fixed, and $\varphi K$ 
being in a small Euclidean neighbourhood of the/some common infinite point $1$ 
of $L_1$ and $L_2$. Then $\varphi K$ does not intersect 
any other $L_i$. However, we may suppose that $K_1$ intersects 
both $L_1$ and $L_2$, and that $\varphi K$ lies on that side of $K_1$, as the
considered 
common infinite endpoint $1$ of $L_1$ and $L_2$. 

Then $(\varphi K) \cap (\psi L)$ is an arc-triangle, with one
infinite vertex. Hence any of its non-trivial symmetries is an axial symmetry,
w.r.t. an axis passing through $1$ (it cannot be a rotation about the
infinite point $1$), and such that $L_1$ and $L_2$, as well as
the base lines of $L_1$ and 
$L_2$, are images of each other by this symmetry. This unique axis of symmetry
must intersect 
the arc-side of $(\varphi K) \cap (\psi L)$ on $K_1$ orthogonally. However, a
small rotation of $K_1$,
about the intersection point of the above axis of symmetry and $K_1$, destroys
this unique symmetry.

{\bf{14.}}
Last we turn to the proof of case 4) from {\bf{8}}.

That is, both graphs consist of vertex-disjoint edges. 
We are going to show that this is impossible. By {\bf{5}}, e.g., 
$L$ has at least
two boundary components.

Let $L_1$ and $L_2$ denote two neighbourly boundary components of $\psi L$,
with $L_2$ following $L_1$ in the positive orientation. 

\newpage

(That is, 
passing on the boundary of $\psi L$, 
taken in one of the models, together with its boundary circle, from $L_1$
to $L_2$, in the positive sense, there are no other connected components of 
bd\,$(\psi L)$, taken in $H^2$, between them.)
Then, denoting by $l_1$ the second 
infinite point of $L_1$, and by
$l_2$ the first infinite point $l_2$ of $L_2$ (both taken in the positive
orientation), the counterclockwise arc
${\widehat{l_1l_2}}$ 
contains no infinite point of any boundary component of $\psi L$.
We may suppose that
the base lines of $L_1$ and $L_2$ are symmetric images of each other w.r.t. the
real axis, with the base line of $L_1$ being in the open lower half-plane, and 
the base line of $L_2$ being in the open 
upper half-plane. 
Let $K_1$ be a boundary component of $\varphi K$. Let 
its infinite end-points be $k_1'$ and $k_1''$, following each other in this
order in the positive sense, on the boundary of $\varphi K$, taken in the
model together with its boundary circle. 
Let us begin with
the position, when $k_1'=l_2$ and $k_1''=l_1$,
and $\varphi K$ lies on the same side of $K_1$, as $1$. 

Now let us translate $\varphi K$, and thus also $K_1$, 
along the real axis a bit, to the left. For the new congruent copy of $K$ 
we will not apply a new notation, but
will preserve the old notation $\varphi K$. Then $k_1'$
and $k_1''$ move a bit (in the conformal model, taken with its boundary
circle).
We want to determine the intersection $(\varphi K) \cap (\psi L)$. 

Let the boundary components of $\varphi K$, or of $\psi L$, be, in the
positive sense, $K_1,K_2,...,K_n$, or $L_1,...,L_m$, respectively. Using the
collinear model, we see that any of $K$ and $L$ can be obtained from a convex
polygon, with all vertices at infinity, whose number of vertices 
is even, by putting hypercycles,
outwards, on each second side, and replacing the remaining sides with the
corresponding arcs of the boundary of the model.
(Including the case when this convex polygon is a $2$-gon, i.e., a segment.)

Then we may suppose that all boundary components of $\varphi K$, except $K_1$ 
(if any), lie strictly on the right hand side of the straight line $l_1l_2$.
All these will be
boundary components of $(\varphi K) \cap (\psi L)$ as well. 
There is still one boundary
component of $(\varphi K) \cap (\psi L)$. This begins at $l_2$, then passes on
$L_2$, then passes on $K_1$, then on some $L_{i(1)}$, then once more on $K_1$,
then on some $L_{i(2)}$, ..., then on some $L_{i(k)}$, then once more on
$K_1$, then on $L_1$, and ends at $l_1$. (One has to observe only that $K_1$
must cross $L_2$, transversally --- observe that both $K_1$ and $L_2$ are
circular arcs in the conformal model --- and then some small arc of it 
still remains in cl\,conv\,$(L_1 \cup L_2)$ a
bit, so that this small arc cannot be "cut off" by any other $L_j$. A similar
reasoning is valid for $L_1$.)

\newpage

Then any non-trivial symmetry of $(\varphi K) \cap (\psi L)$ preserves this
unique non-smooth boundary component of $(\varphi K) \cap (\psi L)$. Thus this
symmetry
is an axial symmetry, which maps $L_1$ to $L_2$, and hence the base line of
$L_1$ to that of $L_2$ --- hence has axis the real 
axis --- and maps the first and last arcs of $K_1$ on 
bd\,$[ (\varphi K) \cap (\psi L) ] $ 
to each other (if there is only one such arc,
then it is mapped to itself). In both cases, the symmetry maps the whole $K_1$
to the whole $K_1$, hence its axis is orthogonal to $K_1$. 

This is no contradiction, since, by construction, $K_1$ is symmetric w.r.t.
the real axis. Now let us consider the point of
intersection 
of $K_1$ with the real axis. Let us rotate a bit $\varphi K$
about this point. Then the combinatorial structure of $(\varphi K) \cap (\psi
L)$ remains of the same type (only possibly the set of indices 
$\{ i(1),...,i(k) \} $ will
change, but this does not invalidate the above considerations). So the unique
non-trivial symmetry has as axis the real axis, that should be 
orthogonal to $K_1$. However, this is
already a contradiction, since, by the above 
rotation, the rotated image
of $K_1$ becomes not orthogonal to the real axis.
$ \blacksquare $
\enddemo


\demo{Proof of Theorem 7} 
{\bf{1.}} We begin with the proof of $(1) \Longrightarrow (2)$ .

{\bf{2.}}
We will make a case distinction.
Either both $K$ and $L$ are strictly convex, or 
one of $K$ and $L$ is not strictly convex. 

We begin with the proof of the case when both $K$ and $L$ are strictly convex.
(Observe that, for $X={\Bbb R}^d$ and $X=H^d$, this follows from the
hypothesis (****) of the theorem.)

{\bf{3.}} 
First, we are going to show that, for any 
$x \in {\text{bd}}\,K$ and any $y \in
{\text{bd}}\, L$, all sectional curvatures exist, and
are equal to some non-negative
constant, and, in case of ${\Bbb R}^d$ and 
$H^d$, even to some positive constant.

{\bf{4.}} Let $n$, $m$ denote the outer unit normals of $K$, or $L$, at 
$x \in {\text{bd}}\,K$,
or $y \in {\text{bd}}\,L$, respectively. 
(Recall that (***) implies smoothness.)
Let us choose an $O \in X$, and let
$e,f$ be opposite
unit vectors in the tangent space of $X$ at $O$. Let us choose $\varphi _0,
\psi _0$,
such that $\varphi _0 x=\psi _0 y=O$, and 
the images (in the tangent bundle) of $n$ or $m$ (by the
maps induced by $\varphi _0$ or $\psi _0$ on the tangent bundle) should be
$e$ or $f$, respectively. Then $(\varphi _0 K) \cap (\psi _0 L) \supset \{ O
\} $. Let $l$ be the geodesic from $O$ in the direction of $e$
(equivalently, of $f$). Let us move $\varphi _0 K , \psi _0 L$ toward each
other, so that their points originally coinciding with $O$ should move on the
straight line $l$, to the respective new positions $O_K$ and $O_L$,
while we allow any rotations of them, about the axis $l$. We denote
these new images by $\varphi K,\,\,\psi L$. 

\newpage

Let the amount of the moving of the points originally coinciding with $O$,  
both for $\varphi _0 K$ and $\psi _0 L$, be a common
distance $OO_K=OO_L=\varepsilon >0$. We may assume that 
$O_K = \varphi x \in $ int$\,(\psi L)$ and 
$O_L = \psi y \in $ int$\,(\varphi K)$.

Then, $C:=(\varphi K) \cap (\psi L)$ has a non-empty interior, and, by strict
convexity of $K$ and $L$, has an
arbitrarily small diameter. Hence it has a centre of symmetry, $c$, say. 
We are going to
show that, for $\varepsilon >0$ sufficiently small, 
$c$ coincides with $O$. First observe that, for $\varepsilon >0$
sufficiently small, we have, by hypothesis (***) of the theorem, 
that the ball of centre $O$ and radius
$\varepsilon $ is contained in $C$. 

{\bf{5.}}
First we deal with the case of $S^d$. Let
$\varphi K'$, or $\psi L'$ denote the half-$S^d$ containing $\varphi K$, or
$\psi L$, and containing
$O_K$, or $O_L$ in its boundary, and thus being there
tangent to ${\text{bd}}\,(\varphi K)$, or ${\text{bd}}\,(\psi L)$, 
respectively. By
$\varphi K \subset \varphi K'$ and $\psi L \subset \psi L'$, we have 
also $C \subset (\varphi K') \cap (\psi L')$. However, 
$(\varphi K') \cap (\psi L')$ contains a unique ball of maximal
radius, namely that with centre $O$, 
and radius $\varepsilon $. Then the same statement holds for $C$ as well.
Thus, the centre
of symmetry $c$ of $C$ must coincide with $O$.
 
Now, we turn to the case of ${\Bbb R}^d$ and
$H^d$. Then, by hypothesis (****) of 
the theorem,
we have that, {\it{in an open $\varepsilon $-neighbourhood of}} 
$x$, {\it{or}} $y$ 
there holds the following
implication. If a point belongs to $K \setminus \{ x \}$, or $L \setminus \{ y
\}$, then it belongs to ${\text{int}}\,K''$ or
${\text{int}}\,L''$, where $K''$ or $L''$ are 
closed balls (for ${\Bbb R}^d$), or closed
convex sets bounded by some hyperspheres (for $H^d$), respectively,
with $x \in $ bd\,$K''$ and $y \in $ bd\,$L''$, and with  bd\,$K''$ and 
bd\,$L''$ having sectional curvatures at most $\varepsilon $.
Moreover, the images of $K''$, or $L''$,
by $\varphi $, or $\psi $, contain
$\varphi (x)=O_K$, or $\psi (y)=O_L$, 
and are there tangent to ${\text{bd}}\,(\varphi K)$, or
${\text{bd}}\, (\psi L)$, and then necessarily
have there their concave sides towards
${\text{int}}\,(\varphi K)$, or ${\text{int}}\,(\psi L)$,
respectively. (Observe, that we may have to decrease $\varepsilon (x)>0$, or
$\varepsilon (y)>0$, from (***) and (****) before Thorem 7, to obtain this.) 

Now we make a case distinction. First we deal with the case $X=H^d$, and second
we will deal with the case $X={\Bbb R}^d$.

So, let $X=H^d$. Without loss of generality, we may
assume, that $K''$ and $L''$ are distance surfaces, with equal distances 
$\varepsilon ' (x) =\varepsilon ' (y) > 0$ 
from their base hyperplanes. Further, we
may suppose $ 0 < \varepsilon < \varepsilon ' (x) =\varepsilon ' (y)$.
We may suppose, that $(\varphi K'') \cap (\psi L'')$ lies in the
intersection of the images by $\varphi $, or $\psi $, of the 
neighbourhoods of $x$, or $y$, mentioned in the beginning of the last but one
paragraph, respectively.
Then, locally, $(\varphi K) \cap (\psi L)$ is
contained in $[{\text{int}}\left( (\varphi K'') \cap (\psi L'') \right) ]
\cup \{ O_K, O_L
\}$. However, then also globally we have the inclusion $(\varphi K) \cap 
(\psi L)
\subset [{\text{int}}\left( (\varphi K'') \cap (\psi L'') \right) ]
\cup \{ O_K, O_L \}$.

\newpage

We show, that the unique ball of maximal radius, contained
in $(\varphi K'') \cap (\psi L'')$, is the one with centre $O$, and radius
$\varepsilon $. Then, the same statement
holds for $(\varphi K) \cap (\psi L)$ as well,
hence the coincidence of $c$ and $O$ will be proved.
 
Observe, that $(\varphi K'') \cap (\psi L'')$ 
is rotationally symmetric about the
axis $O_KO_L$, and is symmetric to the orthogonal halving plane $H$
of the segment $O_KO_L$.
We will show that, if we have a ball, included in $(\varphi K'') 
\cap (\psi L'')$,
with centre different from $O$, and with the centre on the (closed) side of 
the orthogonal halving plane of $O_KO_L$, on which $O_K$ lies, then its radius
is less than $\varepsilon $. (The case of $O_L$ is analogous.) Clearly, we may
restrict ourselves to the case $d=2$.

Let $K'''$ be the base line of $K''$ (i.e., $K''$ is a distance line for
$K'''$). 
Clearly, the straight
line containing $O,O_K,O_L$ is orthogonal to $\varphi K''', H$ and
$\varphi K''$ (these last three curves being distinct, and their
intersections with the straight 
line $O_KO_L$ follow each other in the given order, by
$0< \varepsilon < \varepsilon ' (x)$).
Let the intersection of this straight line with $\varphi K'''$ be $O'$.
The straight lines $\varphi K'''$
and $H$ have no common finite or infinite point. 
The minimal distance of these two straight lines is attained in the
position when we take 
the point $O'$
on $\varphi K'''$, 
and the point
$O$ on $H$. 

Now, let us draw straight lines orthogonal to $\varphi K'''$ 
at each point $O^* \in \varphi K'''$.
Then, for the constant value $\varepsilon ' (x) > 0$, 
we have to pass from any point 
$O^* \in \varphi K'''$ 
a segment of length $\varepsilon ' (x)$ on the respective
orthogonal straight line, towards $\varphi K''$, 
till we reach a point, say, $O^{***}$, on $ \varphi K''$. 
During this motion, we may cross $H$, at some point $O^{**}$. We suppose that
this point $O^{**}$ exists, and, moreover, it lies in $(\varphi K'') \cap (\psi
L'')$.

Then the minimum length of $O^*O^{**}$ is $O'O$, and is attained only
for $O^*=O'$.
Hence, the maximum length of $O^{**}O^{***}$ is $OO_K$, and is attained only
for $O^*=O'$.
Therefore, also the distance of any point $P$, lying on the segment
$O^{**}O^{***}$, to $O^{***}$,
is maximal exactly when 
$O^*=O'$ and $P=O^{**}=O$.
As promised above, this ends the proof, that $c=O$ 
for the case of $H^d$.

There remained the case of $X={\Bbb R}^d$. Then elementary geometrical
considerations
yield that the ball of maximal radius, contained in $(\varphi K'') \cap (\psi
L'')$, has centre $c=O$ (and radius $\varepsilon $). 

{\bf{6.}}
Thus, $(\varphi K) \cap (\psi L)$ 
has as centre of symmetry $O$, and it has a chord
$[O_K,O_L]$, passing through $O$, hence $O_K= \varphi x$ and $O_L= \psi y$ 
are centrally symmetric
images of each other w.r.t. $O$. 
Then the same holds for some of their neighbourhoods, 
relative to ${\text{bd}}\,(\varphi K)$, or ${\text{bd}}\,(\psi L)$,
respectively, for
$\varepsilon $ sufficiently small
(recall that $O_K = \varphi x \in $ int$\,(\psi L)$ and $O_L = \psi y \in $
int$\,(\varphi K)$). 

\newpage

Now, take some $2$-plane containing the straight line
$O_KO_L$. Then, the intersections of some neighbourhoods of $O_K$ and $O_L$,
relative to ${\text{bd}}\,(\varphi K)$, or ${\text{bd}}\,(\psi L)$,
respectively,
with this $2$-plane, are centrally symmetric images of each other. Therefore,
these two curves have, at $O_K$ and $O_L$, the same curvatures 
({\it{sectional curvatures}}), if one of them
exists, or they do not have curvatures there. Observe, that $\varphi $ and
$\psi $ were not determined uniquely, but at their definitions 
there were allowed to apply any
rotations about the axis $l$. Hence, 
either all sectional curvatures (i.e., the 
curvatures of all above curves),
of both $K$
and $L$, at the points $x$ and $y$ are equal, or all of them do not
exist. Observe, 
that $x$ and $y$ were arbitrary points of ${\text{bd}}\,K$ and  
${\text{bd}}\,L$. So either 
\newline 
a) 
all sectional curvatures of both $K$ and $L$ exist, 
at each boundary point of $K$ and $L$, 
and they are equal, namely to some number $\kappa \ge 0$, 
or 
\newline 
b) they do not exist anywhere. 
\newline
However, convex surfaces in ${\Bbb R}^d$
are almost everywhere twice differentiable
(cf. [Sch], pp. 31-32, cited in detail in the sixth paragraph of {\bf{4}} of
the proof of Theorem 2). 
Using the collinear models for $S^d$ and $H^d$, this
holds for $S^d$ and $H^d$ as well.
This rules out possibility b), so possibility a) holds, as
promised above. Clearly, for ${\Bbb R}^d$ and
$H^d$, the hypothesis of our Theorem implies
$\kappa >0$.

{\bf{7.}} Observe, that the above proof also gives, that locally ${\text{bd}}\,
K$, or ${\text{bd}}\,L$ is rotationally symmetric about the normal $n$ at $x$,
or $m$ at $y$, respectively. 
(Recall, that $\varphi $ and $\psi $ were defined only up to
arbitrary rotations about the straight 
line $l$, respectively, and we always had
symmetry about $c=O$.) 
Their $2$-dimensional normal sections,
i.e., the
sections by $2$-planes containing $n$, or $m$, 
are normal sections for all of their points
close to $x$, or $y$, just by local rotational symmetry of ${\text{bd}}\,
K$, or ${\text{bd}}\,L$, respectively. 
Therefore, these $2$-dimensional normal sections have everywhere the same
constant curvature $\kappa $. Hence, locally, these sections 
are congruent cycles in the respective
$2$-dimensional subspaces (for ${\Bbb R}^d$ and
$H^d$ they cannot be straight lines, by 
hypothesis (****) of the
theorem). Therefore, ${\text{bd}}\,K$ and ${\text{bd}}\,L$
are, locally, for $S^d$, congruent 
spheres, including half-$S^d$'s, and, for ${\Bbb R}^d$ and
$H^d$, they are, locally, congruent
spheres, paraspheres, or congruent hyperspheres (they cannot be hyperplanes). 
Thus, locally, any of ${\text{bd}}\,K$ and 
${\text{bd}}\,L$ is an analytic surface, given up to congruence. 

Now, let $x \in {\text{bd}}\,K$ be arbitrary. For some relatively
open geodesic $(d-1)$-ball $B_x$
on ${\text{bd}}\,K$, with centre $x$, we have that $B_x$
is a subset of an above analytic hypersurface; if the above hypersurfaces are
spheres, then we assume that the $B_x$'s are
at most half-spheres of these spheres. 

\newpage

For $x_1,x_2 \in {\text{bd}}\,K$, with $
B_{x_1} \cap
B_{x_2} \ne \emptyset $, we have that $B_{x_1}$ and 
$B_{x_2}$ are subsets of the same analytic hypersurface, i.e.,
they are open 
subsets of the same sphere, parasphere, or hypersphere. Now, let us
introduce an equivalence relation on the points $x$ of ${\text{bd}}\,K$. Two
such points $x',x''$ are called {\it{equivalent}}, 
if there exists a finite sequence
$x'=x_1, \ldots , x_n=x'' \in {\text{bd}}\,K$, such that $
B_{x_i} \cap B_{x_{i+1}} \ne \emptyset $, for each
$i=1,\ldots , n-1$. Clearly, the union of each equivalence class is a
relatively open
subset of a sphere, parasphere or hypersphere, and also is relatively open in 
${\text{bd}}\,K$. Thus, they form a relatively 
open partition of ${\text{bd}}\,K$, 
which implies, that they
form a relatively 
open-and-closed partition of ${\text{bd}}\,K$. Thus, the union of each
equivalence class is the union of some components of ${\text{bd}}\,K$. Clearly,
no $B_x$ can intersect different connected components of ${\text{bd}}\,K$,
since $B_x$ is connected. Hence, the unions of the
equivalence classes are subsets of some 
connected components of ${\text{bd}}\,K$. Since also they are unions of some
connected components of ${\text{bd}}\,K$, they are exactly the connected
components of ${\text{bd}}\,K$.

Up to now, we know the following. The connected components of
${\text{bd}}\,K$, and of ${\text{bd}}\,L$, are relatively open subsets of some
congruent spheres/pa\-ra\-spheres/hy\-per\-spheres. 
Since ${\text{bd}}\,K$ is closed
in $X$, its connected components, being relatively closed in ${\text{bd}}\,K$,
are closed in $X$ as well. Thus, the connected components of ${\text{bd}}\,K$
are non-empty, relatively open-and-closed subsets of some congruent 
spheres/pa\-ra\-spheres/hy\-per\-spheres. However, spheres, paraspheres and
hyperspheres are connected, i.e., have no non-empty, 
relatively open-and-closed proper 
subsets. Therefore, the connected components of
${\text{bd}}\,K$, and, similarly, of ${\text{bd}}\,L$, are congruent 
spheres/pa\-ra\-spheres/hy\-per\-spheres.

This shows the implication $(1) \Longrightarrow (2)$ for the case, when both
$K$ and $L$ are strictly convex, i.e., in the first case in {\bf{2}}.

{\bf{8.}}
Now suppose that 
one of $K$ and $L$ 
is not strictly convex, that is the second case in {\bf{2}}. By hypothesis
(****) of the theorem, this can happen only for $X=S^d$. 

By (***) both $K$ and $L$ are smooth. We consider two cases for $K$ (and
analogously for $L$). We have either diam$\,K < \pi $, or diam$\,K= \pi $.

In the first case, $K$ is contained in an open half-sphere. Let us suppose that
this half-sphere is the southern half-sphere. Then the collinear model is
defined in a neighbourhood of $K$, and the image of $K$ is a compact convex
set in the model ${\Bbb R}^d$. Such a set has an exposed point, i.e., a point
$z$ such that $\{ z \} $ 
is the intersection of the image of $K$ and a hyperplane in the model
${\Bbb R}^d$, cf. [Sch], Theorem 1.4.7 (Straszewicz's theorem). 

\newpage

In the second case, $K$ contains two antipodal points of $S^d$, and
we may suppose
that these are $(0,...,0, \pm 1)$. Since $K$ is smooth at $(0,...,0,1)$,
therefore we may suppose that it has at $(0,...,0,1)$
the tangent hyperplane (in $S^d$)
$\{ (x_1,...,x_d,x_{d+1}) \in S^d \mid x_1 = 0 \} $, and $K$ lies on the side 
$\{ (x_1,...,x_d,x_{d+1}) \in S^d \mid x_1 \ge 0 \} $ of this hyperplane.
Clearly $K$ consists of entire half-meridians, connecting $(0,...,0, \pm 1)$.
By the hypothesis about the tangent hyperplane, each half-meridian, whose
relative interior lies in the
open half-sphere, given by $x_1 > 0 $, lies entirely in $K$. 
Therefore, $K$ contains the closed half-sphere, given by $x_1 \ge 0$. Since, by
hypothesis (**) of the theorem, we have 
$K \ne S^d$, we have that $K$ is a half-sphere.

Considering also $L$, we have also
that either $L$ has an exposed point, or $L$ is
a half-sphere. So, unless both $K$ and $L$ are half-spheres --- when 
we are done --- we have that, 
e.g., $K$ has an exposed point $x$. Then let $y \in $ bd$\,L$. Now
we can repeat the procedure described in {\bf{4}}. Then $(\varphi K) \cap
(\psi L)$ has an arbitrarily 
small diameter, hence is centrally symmetric by (1).
Then we have the situation decribed in {\bf{4}} and the first paragraph of
{\bf{5}}. Then the
first two sentences of {\bf{6}} are valid also here. 
That is, some small neighbourhoods of $O_K=\varphi x$, or $O_L=\psi y$,
relative to bd$\,(\varphi K)$, or to bd$\,(\psi L)$, respectively, are
centrally symmetric images of each other w.r.t. $O$, with 
$O_K=\varphi x$ and $O_L=\psi y$ being the 
centrally symmetric images of each other w.r.t. $O$.

This implies, that also $y$ is an exposed point of $L$ (observe that to be an
exposed point of a closed convex set is a local property). 
That is, all boundary
points of $L$ are exposed points of $L$. 
Now, changing the roles of $K$ and $L$, we
obtain, that also all boundary points of $K$ are exposed points of $K$.
In other words, both $K$ and $L$ are strictly convex.
However, this contradicts the hypothesis in the first sentence of {\bf{8}}.

{\bf{9.}}
The implication $(2) \Longrightarrow (1)$ of the theorem is proved by copying
the respective proof from {\bf{7}} of the proof of Theorem 2. Now we will
have central
symmetry, since in (2) of the theorem we have congruent connected components.

In fact, the intersection of two congruent balls (with non-empty interior)
is centrally symmetric.

A compact intersection of two paraballs $\varphi K$ and $\psi L$ 
(with non-empty interior)
is centrally symmetric. In fact, the
infinite points of the two paraballs, say, $k$ and $l$, are different. 
We consider the straight
line $kl$. Let the other points of bd\,$(\varphi K)$ and bd\,$(\psi L)$ 
on $kl$ be $k'$ and $l'$.

\newpage

We may suppose that the order of the points on $kl$ is $k,l',k',k$. Then the
midpoint of the segment $k'l'$ is the centre of symmetry of $(\varphi K) \cap
(\psi L)$.

For the case when the boundary components are congruent hyperspheres, these
hyperspheres are distance surfaces for some distance $c>0$. Like in the 
proof of {\bf{7}} of the proof of Theorem 2, we may restrict ourselves to the
case when $(\varphi K) \cap (\psi L)$ is bounded only by one boundary
component of $\varphi K$, and of $\psi L$. (Even the different boundary
components have here distances at least $2c$.) That is, we have a compact
intersection (with non-empty interior)
of two convex sets, $\varphi K$ and $\psi L$, 
bounded by congruent hypersperes. Then the sets of infinite points of 
$\varphi K$ and $\psi L$ are disjoint.
Considering the collinear model, this implies that
the base hyperplanes of bd\,$(\varphi K)$ and bd\,$(\psi L)$ 
have no finite, or infinite
points in common. Let us consider the segment realizing the distance of these
hyperplanes. Then its 
midpoint is the centre of symmetry of $(\varphi K) \cap (\psi L)$.
$ \blacksquare $
\enddemo


\demo{Proof of Theorem 8}
{\bf{1.}}
We have to prove only $(1) \Longrightarrow (2)$.
 
Observe that (1) of Theorem 8 implies (1) of Theorem 7, and (1) of Theorem 7
implies, by Theorem 7, that the connected components of the boundaries both of
$K$ and $L$ are 
either 1) congruent 
spheres (for $X=S^d$ of radius at most $\pi /2$),
or 2) paraspheres, or 3) congruent hyperspheres.

In case 1) $K$ and $L$ are congruent balls, hence $(2)$ is proved.

There remained the cases when we have 
$X=H^d$, and $K$ and $L$ are 2) two
paraballs, or 3) the boundary components both of $K$ and $L$ 
are congruent
hyperspheres, and their numbers are at least $1$, but at most countably
infinite.
We will copy the respective parts of the proof of Theorem 5, {\bf{3}} 
and {\bf{5}}. We are going to show, that neither of these cases can occur.

In case 2) $K$ and $L$ are paraballs. We choose
$\varphi $ and $\psi $ so, that $\varphi K = \psi L$. Then their intersection
is a paraball, that is not centrally symmetric, like at Theorem 5, {\bf{3}}.

We turn to case 3).
It will be convenient to use the conformal model. 
Let all boundary
components $K_i$ of $\varphi K$, and $L_i$ of $\psi L$, 
be congruent hyperspheres, with base 
hyperplanes $K_{0,i}$ and $L_{0,i}$. 
Denote by
$l$ the common value of the distance, for which these hyperspheres are
distance surfaces for their base hyperplanes. 
(By hypothesis (****) of the theorem
we have $l>0$.)
These base hyperplanes 
bound closed convex sets $K_0$, or $L_0$, possibly with empty interior,
not containing any of the hyperspheres $K_i$, or $L_i$, 

\newpage

and such that 
the parallel domain of $K_0$, or $L_0$, with distance $l$,
equals $\varphi K$, or $\psi L$, respectively. Cf.
the proof of Theorem 5, {\bf{5}}. 

We choose such positions of $\varphi K$ and
$\psi L$, 
that $K_{0,1}$ and $L_{0,1}$ project to a copy of $H^2$ in $H^d$ (with a
projection along straight lines orthogonal to the copy of $H^2$), so that
their projections are like $K_{0,1}$ and $L_{0,1}$ from the proof of Theorem
5, {\bf{5}}. For simplicity, here we assume that this copy of $H^2$ contains
the centre of the model, and the axis of symmetry of the above
projections, in this copy of $H^2$, 
passes through the centre of the model. (This implies that, in the
conformal model, these projections have the same lengths.) 
Then, the proof of Theorem 5, {\bf{5}} gives, that $(\varphi K) \cap (\psi L)$ 
is the intersection of two closed convex sets, bounded by the congruent
hyperspheres $K_1$ and $L_1$. 

Now, following the proof
of Theorem 5, {\bf{5}}, we will show that
this intersection is not centrally symmetric.
In fact,
in the conformal model, the hyperspheres 
$K_1$ and $L_1$ are subsets
of spherical surfaces (congruent, in the Euclidean sense, 
in the conformal model), 
with their centres $k$ and $l$ in the
Euclidean plane spanned by the above (conformal) model circle of $H^2$.
Moreover, in the conformal model,
their intersection is a $(d-2)$-sphere, of (Euclidean) radius less than $1$,
that is (in Euclidean sense) rotationally
symmetric about the straight line $kl$, and
touches the boundary of the model ball
at one point (namely at 
the intersection of the projections of $K_{0,1}$ and $L_{0,1}$), 
and has all other points 
in the model. Hence, in the conformal model, 
the intersection of the closed convex sets bounded by
$K_1$ and $L_1$ also is (in Euclidean sense)
rotationally symmetric about the axis $kl$, and 
also touches the boundary of the model at one point, and has all other points 
in the model (in fact, it is contained in the Thales $(d-1)$-sphere 
of the above $(d-2)$-sphere, taken in the model). 
Therefore this intersection cannot be centrally symmetric. 
$ \blacksquare $
\enddemo


\definition{Acknowledgements} 
The authors express their gratitude to I. B\'ar\'any, for
carrying the problem, and bringing the two authors
together; to V. Soltan, for having sent to the second named
author the manuscript of [So2], 
prior to its publication; and to K. B\"or\"oczky
(Sr.), for calling the attention of the second named author 
to the fact that, in Theorem
2, our original hypothesis $C^3_+$ was unnecessary. Following K. B\"or\"oczky's
arguments, the
authors finally succeded to eliminate all smoothness hypotheses. 
W
\enddefinition




\newpage

\Refs

\widestnumber\key{WWW}


\ref 
\key 
\book 
\by  
\publ 
\publaddr 
\yr 
\endref 

\ref 
\key  
\by 
\paper  
\jour 
\pages  
\endref   

\ref
\key 
\by 
\paper 
\jour 
\vol 
\yr 
\pages  
\endref 

\ref 
\key AVS 
\book Geometry of spaces of constant curvature
\by D. V. Alekseevskij, E. B. Vinberg, A. S. Solodovnikov
\publ Geometry II (Ed. E. B. Vinberg), Enc. Math. Sci. 29, 1-138, Springer, 
\publaddr Berlin
\yr 1993
\MR {\bf{95b:}}{\rm{53042}}
\endref 

\ref  
\key Ba
\book Nichteuklidische Geometrie, Hyperbolische Geometrie der Ebene 
\by R. Baldus
\publ 4-te Aufl., Bearb. und erg\"anzt von F. L\"obell, Sammlung G\"oschen
970/970a, de Gruyter
\publaddr Berlin 
\yr 1964
\MR {\bf{29\#}}{\rm{3936}}
\endref 

\ref 
\key BF 
\book Theorie der konvexen K\"orper, {\rm{Berichtigter Reprint}}
\by T. Bonnesen, W. Fenchel
\publ Springer
\publaddr Berlin-New York
\yr 1974
\MR {\bf{49\#}}{\rm{9736}}
\endref 

\ref 
\key Bo
\book Non-Euclidean geometry, a critical and
historical study of its developments, {\rm{With a Supplement containing the
G. B. Halstead translations of}} ``The science of absolute space'' {\rm{by 
J. Bolyai and}} ``The theory of parallels'' {\rm{by N. Lobachevski}}
\by R. Bonola
\publ Dover Publs. Inc.
\publaddr New York, N.Y.
\yr 1955
\MR {\bf{16-}}{\rm{1145}}
\endref 

\ref 
\key C
\book Non-Euclidean Geometry, 6-th ed. 
\by H. S. M. Coxeter 
\publ Spectrum Series, The Math. Ass. of America
\publaddr Washington, DC
\yr 1998
\MR {\bf{99c:}}{\rm{51002}}
\endref 

\ref 
\key HM
\by E. Heil, H. Martini
\paper Special convex bodies 
\jour 
In: Handbook of Convex
Geometry, (eds. P. M. Gruber, J. M. Wills), North-Holland,
Amsterdam etc., 1993, Ch. 1.11
\pages 347-385 
\MR {\bf{94h:}}{\rm{52001}}
\endref   
 
\ref
\key H
\by R. High
\paper Characterization of a disc, Solution to problem 1360 (posed by
P. R. Scott)
\jour Math. Magazine
\vol 64 
\yr 1991
\pages 353-354 
\endref 

\ref 
\key J-CM
\by J. Jer\'onimo-Castro, E. Makai, Jr.
\paper 
\jour in preparation 
\vol 
\yr 
\pages  
\endref 

\ref 
\key L
\book Nichteuklidische Geometrie, {\rm{3-te Auflage}}
\by H. Liebmann
\publ de Gruyter
\publaddr Berlin
\yr 1923
Jahr\-buch Fortschr. Math. {\bf{49}}, 390
\endref 

\ref 
\key P
\book Nichteuklidische Elementargeometrie der Ebene, {\rm{Math. Leitf\"aden}}
\by O. Perron 
\publ Teubner
\publaddr Stuttgart
\yr 1962
\MR {\bf{25\#}}{\rm{2489}}
\endref 

\ref 
\key Sch
\book Convex bodies: the Brunn-Minkowski theory, {\rm{Encyclopedia of}}
\newline
{\rm{Math. and its Appls., Vol.}} {\bf{44}}
\by R. Schneider 
\publ Cambridge Univ. Press
\publaddr Cambridge
\yr 1993
\MR {\bf{94d:}}{\rm{52007}}
\endref 

\ref
\key So1
\by V. Soltan
\paper Pairs of convex bodies with
centrally symmetric intersections of translates
\jour Discrete Comput. Geom. 
\vol 33
\yr 2005
\pages 605-616 
\MR {\bf{2005k:}}{\rm{52012}} 
\endref 

\ref
\key So2
\by V. Soltan
\paper Line-free convex bodies
with centrally symmetric intersections of translates 
\jour Revue Roumaine de Math. Pures et Appl. 
\vol 51
\yr 2006
\pages 111-123 
\MR {\bf{2007k:}}
\newline
{\rm{52010}}.
{\rm{Also in: Papers on Convexity and Discrete geometry, Ded. to T. 
Zamfirescu on the occasion of his
60-th birthday, Editura Academiei Rom\^ane, Bucure\c sti, 2006, 411-423}}
\endref 

\ref 
\key St
\book Differential Geometry
\by J. J. Stoker 
\publ New York Univ., Inst. Math. Sci.
\publaddr New York
\yr 1956
\endref 

\ref
\key V
\by I. Vermes
\paper \"Uber die synthetische Behandlung der
Kr\"ummung und des Schmieg\-zy\-kels der ebenen Kurven in der
Bolyai-Lobatschefskyschen Geometrie
\jour Stud. Sci. Math. Hungar.
\vol 28
\yr 1993
\pages 289-297 
\MR {\bf{95e:}}{\rm{51030}} 
\endref 

\endRefs
\enddocument